\RequirePackage{fix-cm}
\documentclass[smallextended]{svjour3}

\usepackage{cmap} 

\smartqed  

\usepackage[utf8]{inputenc}

\usepackage{color}
\usepackage{soul} 
\usepackage{bm} 

\usepackage{paralist}

\usepackage{graphicx}
\graphicspath{ {./img/} }
\DeclareGraphicsExtensions{.pdf}
\usepackage{subcaption}

\usepackage{amsmath}
\usepackage{amssymb}
\usepackage{mathtools}
\usepackage{amsrefs}

\newcommand{\abs}[1]{ \ensuremath{\left\lvert{#1}\right\rvert}}

\newcommand{\conj}[1]{ \ensuremath{\overline{#1}}}

\renewcommand{\theta}{\vartheta}

\usepackage[shortcuts]{extdash}
\DeclareMathOperator{\sgn}{sgn}
\DeclareMathOperator{\Log}{Log}
\DeclareMathOperator{\Arg}{Arg}
\DeclareMathOperator*{\ilim}{i-lim}
\DeclareMathOperator*{\elim}{e-lim}

\DeclareMathOperator*{\argmax}{arg\,max}

\DeclareMathOperator{\Arccot}{Arccot}

\def\Xint#1{\mathchoice
   {\XXint\displaystyle\textstyle{#1}}%
   {\XXint\textstyle\scriptstyle{#1}}%
   {\XXint\scriptstyle\scriptscriptstyle{#1}}%
   {\XXint\scriptscriptstyle\scriptscriptstyle{#1}}%
   \!\int}
\def\XXint#1#2#3{{\setbox0=\hbox{$#1{#2#3}{\int}$}
     \vcenter{\hbox{$#2#3$}}\kern-.5\wd0}}

\def\dashint{\Xint-}

\usepackage{hyperref}
\hypersetup{
bookmarksopen=True,
bookmarksnumbered=True,
breaklinks=True
}

\usepackage{tikz}
\usetikzlibrary{matrix}
\usetikzlibrary{arrows}
\usetikzlibrary{positioning}
\usepackage{relsize}
\tikzset{fontscale/.style = {font=\relsize{#1}}}

\begin{document}

\title{Conditioning moments of singular measures for entropy maximization II: Numerical examples \thanks{To Ed Saff on the occasion of his seventieth birthday.}}
\titlerunning{Conditioning moments for MAXENT}

\author{Marko Budi\v{s}i\'{c} \and Mihai Putinar}
\authorrunning{M. Budi\v{s}i\'{c} \and M. Putinar} 
\institute{Department of Mathematics,
  University of Wisconsin,
  Madison, WI, USA
  \email{marko@math.wisc.edu} \and
Department of Mathematics,
University of Santa Barbara,
Santa Barbara, CA, USA \email{mputinar@math.ucsb.edu}, and
School of Mathematics \& Statistics
Newcastle University, Newcastle upon Tyne
NE1 7RU, United Kingdom.
\email{Mihai.Putinar@newcastle.ac.uk}}

\date{March 24, 2015}
\maketitle

\begin{abstract}
  If moments of singular measures are passed as inputs to the entropy
  maximization procedure, the optimization algorithm might not terminate. The
  framework developed in~\cite{budisic_conditioning_2012} demonstrated how
  input moments of measures, on a broad range of domains, can be conditioned
  to ensure convergence of the entropy maximization. Here we numerically
  illustrate the developed framework on simplest possible examples: measures
  with one-dimensional, bounded supports. Three examples of measures are used
  to numerically compare approximations obtained through entropy maximization
  with and without the conditioning step.
  \\
  MSC 30E05,30E20,41A46.  
\end{abstract}

\section{Introduction}
\label{sec:introduction}

Among all inverse problems, old and new, theoretical or applied, the
reconstruction or, at least, approximation of a positive measure from moment
data stands as a major topic of continued interest with history of a century
and a half. Modern motivations for studying moment problems arise in continuum
mechanics, statistics, image processing, control theory, geophysics,
polynomial optimization, to name only a few. Added to these are purely
theoretical aspects of convex algebraic geometry or harmonic analysis. While
most cases deal either with measures absolutely continuous with respect to
Lebesgue measure or with purely atomic measures, the class of continuous
singular measures has been poorly studied from the point of view of moment
data inversion.  The present note is a continuation of
article~\cite{budisic_conditioning_2012} which developed a theoretical
framework that includes all singular measures in the ubiquitous maximum
entropy closure. Our aim here is to numerically demonstrate on several simple
examples how input moment data can be conditioned using a nonlinear triangular
transform, with the purpose of stimulating research into effects of such
moment conditioning on the numerical accuracy, algorithm implementation, and
even theoretical underpinnings of the maximum entropy method.  For
illustration, we choose the comfort of the measures with one-dimensional,
bounded supports, although our previous work~\cite{budisic_conditioning_2012}
additionally covered unbounded supports and multivariate moments. A
forthcoming third article in this series will deal with inverse moment
problems involving singular measures appearing in analysis of spectra of some
concrete dynamical systems~\cite{budisic_applied_2012}.

The amount of input moment data available is one of the distinguishing
characteristics between theoretical and applied problems. Applied problems
almost exclusively feature limited input data (\emph{the truncated moment
  problem}) which means that the inversion is under-determined: there will be
many measures whose moments match the inputs. The unique solution measure is
obtained through \emph{moment closure}, i.e., by posing an additional
assumption on the input measure, which constrains the set of feasible measures
in which the solution is searched for and ensures a unique solution. Requiring
that the solution maximizes Boltzmann-Shannon entropy
(MAXENT)\footnote{Depending on the convention for choosing the sign of the
  entropy functional, both ``minimization'' and ``maximization'' appear in
  literature.} is an example of a closure that is particularly common in
physical sciences and statistics~\cite{borwein_maximum_2012}.

Explanations of intuition behind MAXENT approach abound in literature, ranging
from rigorous to philosophical~\cites{jaynes_rationale_1982, csiszar_why_1991,
  grendar_maximum_2001, seidenfeld_entropy_1986}. Informally, MAXENT
probability measures are those that match the input moment data, but are
otherwise unbiased towards any particular outcome. For our purposes, the
important feature of MAXENT measures is that they are continuous, i.e., they
have a density function.  Their densities are in the form of exponential
polynomials \(e^{a_{n}x^{n} + a_{n-1}x^{n-1}+\dots}\) which means that the
``rougher'' the measure that generated the input data, the less likely it is
for MAXENT to provide a good approximation from finitely-many moments. In the
extreme, when input moments come from a singular measure, MAXENT approach may
even fail to converge~\cite{budisic_conditioning_2012}. Singular measures,
however, are hardly pathological if a broad range of domains is surveyed. In
particular, they often appear as invariant measures of dynamical systems whose
dynamics settles on a lower-dimensional attractor. In function theory they
naturally arise in the study of boundary values of bounded analytic functions
in a planar region~\cite{goluzin_geometric_1969}.

In Section~\ref{sec:analytics}, we briefly re-state the basic theory
from~\cite{budisic_conditioning_2012} needed to understand the
single-variable, bounded-domain problem. Section~\ref{sec:numerics} discusses
three numerical problems implementation of the algorithm in
MATLAB~\cite{the_mathworks_inc_matlab_2014} and \texttt{Chebfun}
~\cite{driscoll_Chebfun_2014}.  Finally, we summarize the results in
Section~\ref{sec:conclusions}.

\section{Analytical description}
\label{sec:analytics}

\subsection{Three representations of a measure}
\label{sec:three-representations}
Consider the bounded interval \([-\pi, \pi]\) and interpret it as the boundary
\(\partial \mathbb{D}\) of the open unit disk \(\mathbb{D}\).  Endow the
interval with a measure \(\mu\) with no assumptions on regularity of \(\mu\)
and, in particular, on existence of density \(\mu'(\theta)\).  The complex
trigonometric moments of \(\mu\), also known as Fourier coefficients, are
given by
\begin{equation}
  \label{eq:trig-moments}
  \tau_{\mu}(k) = \frac{1}{2\pi} \int_{-\pi}^{\pi} e^{-i k \theta}d\mu(\theta).
\end{equation}
If the real-valued density \(\mu'(\theta)\) is square-integrable,
trigonometric moments are featured in its Fourier expansion:
\begin{equation}
\label{eq:trig-moments-fourier-coeffs}
\begin{aligned}
  \mu'(\theta) &= \sum_{n = -\infty}^{\infty} \tau_{\mu}(n)e^{in\theta} \\&=
  \tau_{\mu}(0) + 2\sum_{n=1}^{\infty}[ \Re \tau_{\mu}(n) \cos(n\theta) -
  \Im \tau_{\mu}(n)\sin(n\theta)].
\end{aligned}
\end{equation}

The measure \(\mu\) induces the function \(\mathring \mu(z) : \mathbb{C} \to \mathbb{C}\) through the Cauchy integral\footnote{In literature, it is common to see \(2\pi i\) in the denominator, which makes formulas in this paper slightly unorthodox.}
\begin{equation}
  \label{eq:analytic-representation}
  \mathring{\mu}(z) := \frac{i}{2\pi} \int_{-\pi}^{\pi} \frac{d \mu(\theta)}{1 - e^{-i\theta}z} = \frac{1}{2\pi} \int_{\zeta \in \partial \mathbb D} \frac{d\mu(\zeta)}{\zeta - z}.
\end{equation}
The Cauchy transform is commonly defined on the real line as an integral of the Cauchy kernel \({(x-z)}^{-1}\) against a measure. As the expression~\eqref{eq:analytic-representation} is equivalent to the classical Cauchy transform on the real line, we refer to it as the Cauchy transform as well. The function \(\mathring{\mu}(z)\) is known as the \emph{analytic representation} of \(\mu\)~\cite[\S 10.8]{king_hilbert_2009}, as \(\mathring{\mu}\) is analytic everywhere except on the unit circle \( \partial \mathbb{D}\). Its derivatives at \(z = 0\) are given by
\begin{equation}
  \label{eq:taylor-coefficients}
  \frac{d^{k} \mathring{\mu}}{dz^{k}} = \frac{i k!}{2\pi}\int_{-\pi}^{\pi} \frac{e^{-ik\theta} d\mu(\theta)}{(1 - e^{-i\theta}z)^{k+1}}, \quad
  \frac{d^{k} \mathring{\mu}}{dz^{k}}(0) = i\, k!\, \tau_{\mu}(k),
\end{equation}
so it follows that the power expansion of the analytic representation \(\mathring{\mu}(z)\) on the unit disk \(\mathbb{D}\) encodes the trigonometric moment sequence \(\mu\):
\begin{equation}
  \label{eq:power-expansion}
  \mathring{\mu}(z) = i \sum_{k = 0}^{\infty} \tau_{\mu}(k) z^{k}.
\end{equation}

When \(\mu(\theta)\) has a H\"older-continuous density \(\mu'(\theta)\), it is possible to recover it from \(\mathring \mu(z)\) by taking limits of \(z \to e^{i\theta}\) along directions not tangential to the unit circle. Depending on the region approach with respect to the unit disk, we distinguish interior and exterior non-tangential limits, respectively denoted by \(\ilim_{z \to e^{i\theta}} \mathring{\mu}(z) \) and \(\elim_{z \to e^{i\theta}} \mathring{\mu}(z)\). For H\"older-continuous \(\mu'(\theta)\), the non-tangential limits depend only on the region (interior/exterior) of the approach, but not on the precise direction. While taking \(\ilim\) is just a standard limit along the appropriate radius of \(\mathbb D\), expressions~\eqref{eq:analytic-representation} and~\eqref{eq:power-expansion} cannot be used directly to evaluate \(\elim_{z \to e^{i\theta}}\) where approach is through the complement of \(\overline{\mathbb{D}}\). However, it is possible to define the \emph{associate function} \(\mathring{\mu}^{\dagger}(z) := \conj{\mathring{\mu}(1/\conj{z})}\), analytic in \(\mathbb D\), whose \(\ilim\) limit relates to the \(\elim\) of \(\mathring{\mu}(z).\) Full details are available in~\cite[\S 14.2.I]{henrici_applied_1977}, but for our purposes the following formula for \(\elim\) will suffice:
\begin{equation}
  \label{eq:elim-as-ilim}
  \elim_{z \to e^{i\theta}} \mathring{\mu}(z) = i \tau_{\mu}(0) + \conj{\ilim_{z \to e^{i\theta}} \mathring{\mu}(z)}.
\end{equation}

Pointwise values of limits are given by Plemelj--Sokhotski formulas\footnote{These formulas are commonly labeled by various subsets of names of Plemelj, Sokhotski, and Privalov, all of which have completed substantial work on them. Typically, these formulas are stated in a complex-argument form. The crucial step in conversion to angle form on the unit circle is the identity \((1-e^{ix})^{-1} = [1 + i \cot (x/2)]/2\).}~\cite[\S 14.2-IV]{henrici_applied_1977}:
\begin{equation}
  \label{eq:plemelj-sokhotski-formulas}
\begin{aligned}
  \ilim_{z \to e^{i\theta}} \mathring{\mu}(z) &= \phantom{-}\frac{i}{2} \mu'(\theta) + \frac{i}{2}\tau_{\mu}(0) - \frac{1}{2}\mathcal{H} \mu'(\theta) \\
  \elim_{z \to e^{i\theta}} \mathring{\mu}(z) &= -\frac{i}{2} \mu'(\theta) + \frac{i}{2}\tau_{\mu}(0) - \frac{1}{2}\mathcal{H} \mu'(\theta)
\end{aligned}
\end{equation}
where \(\mathcal H \mu'\) is the Hilbert transform of the density \(\mu'(\theta)\)
\begin{equation}
  \label{eq:hilbert-transform}
  \mathcal H \mu'(\theta) = \frac{1}{2\pi} \dashint_{-\pi}^{\pi} \cot\frac{\theta-\alpha}{2}\mu'(\alpha) d\alpha
\end{equation}
evaluated using the principal value integral \(\dashint\).
Subtracting the formulas we recover the density \(\mu'(\theta)\) pointwise:
\begin{equation}
  \label{eq:plemelj-sokhotski-density}
  \begin{aligned}
    \mu'(\theta) &= \phantom{-}i \left[\elim_{z \to e^{i\theta}}
      \mathring{\mu}(z) - \ilim_{z \to e^{i\theta}}
      \mathring{\mu}(z)\right] \\
    \mathcal{H} \mu'(\theta) &= \phantom{i}-\left[\elim_{z \to e^{i\theta}}
      \mathring{\mu}(z) + \ilim_{z \to e^{i\theta}}
      \mathring{\mu}(z)\right] +i\tau_{\mu}(0).
  \end{aligned}
  \end{equation}

\begin{figure}[htb]
  \centering
  \begin{tikzpicture}[
    description/.style={fill=white, inner sep=2pt},
    mathnode/.style={execute at begin node=\(\scriptstyle, execute at end node=\)},
    textnode/.style={execute at begin node=\footnotesize},
    ge/.style={fontscale=1.5}]
\matrix (m)  [matrix of math nodes, row sep=4em,
column sep=6em, minimum width=2em, cells=ge, anchor=center] 
{ \mu(\theta) &                    \\
      &  \mathring{\mu}(z) \\
  \tau_{\mu}(k) & \\};
\path[-]
(m-1-1) edge [thick,->,transform canvas={yshift=0.75ex}] node [sloped, anchor=center,above,textnode] {Cauchy transform}
(m-2-2)
(m-2-2) edge [dashed,thick,->,transform canvas={yshift=-.75ex}] node [sloped, anchor=center,below,textnode] {Plemelj--Sokhotski}
(m-1-1)
(m-1-1) edge [thick,->,transform canvas={xshift=-.75ex}] node [left,mathnode] {\int e^{-ik\theta}d\mu}
(m-3-1)
(m-3-1) edge [dotted,thick,->,transform canvas={xshift=0.75ex}] node [right,textnode] {\parbox[c]{4em}{Inverse problem}}
(m-1-1)
(m-2-2) edge [thick,->,transform canvas={yshift=-.75ex}] node [right, sloped, anchor=center,below,textnode] {Power expansion}
(m-3-1)
(m-3-1) edge [thick,->,transform canvas={yshift=0.75ex}] node [sloped, anchor=center,above,textnode] {Formal series}
(m-2-2)
;
  \end{tikzpicture}
  \caption{Relations between the measure \(\mu(\theta)\), its analytic representation \(\mathring\mu(z)\), and trigonometric moments \(\tau_{\mu}(k)\). Plemelj--Sokhotski formulas (dashed) apply only when \(\mu\) has a H\"older-continuous density function. A unique solution to the inverse problem (dotted) is available in certain cases when the full sequence \(\tau_{\mu}(k),\ \forall k \in \mathbb{Z}\) is known or, by closure, under additional restrictions on \(\mu\).}
  \label{fig:three-representations}
\end{figure}
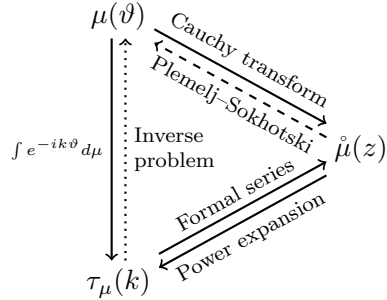

\autoref{fig:three-representations} illustrates how the concepts of measure \(\mu\), its moment sequence \(\tau_{\mu}(k)\), and analytic representation \(\mathring{\mu}(z)\) relate. The stated relationships are well-documented in classical literature~\cites{henrici_applied_1977,goluzin_geometric_1969,king_hilbert_2009} and further references listed in~\cite{budisic_conditioning_2012}.

As mentioned in Section~\ref{sec:introduction}, our goal is to specify a closure for the truncated moment problem: recover a representation of \(\mu\) from a partial knowledge of moments \(\tau_{\mu}\) or, equivalently, truncation of power series of \(\mathring{\mu}(z)\). Two well-known examples of closures are the Pad\'e approximation and the entropy maximization, both requiring additional assumptions on the (non-)existence of density \(\mu'(\theta)\) and its smoothness. The Pad\'e approximation can treat both singular and continuous measures \(\mu\), but at the cost that it always gives an atomic measure as the result which is, in practice, a poor approximation of any measure with a continuous density. On the other hand, the maximum entropy closure does not converge when applied to moments generated by some singular measures, e.g., a point mass measure~\cite{budisic_conditioning_2012}, thus not producing an approximation without an additional stopping criterion.

\subsection{Even singular measures are represented by bounded densities}
\label{sec:singular-measures-by-densities}

In this section we relate \(\mathring{\mu}\) to another analytic function with very similar properties but with an important improvement: it is represented by a Cauchy integral of a measure with a \emph{bounded} density, even in cases of singular \(\mu\). We follow ideas of Markov~\cite{akhiezer_questions_1962} and Aronszajn and Donoghue~\cite{aronszajn_exponential_1956} in this development.

For positive measures \(\mu\), \(\mathring{\mu}(z)\) maps \(\mathbb D\) into the upper half plane \(\mathbb H^{+}\), which is easy to demonstrate by computing the imaginary part of the integral kernel \(i/(1 - e^{-i\theta}z)\). The consequence is that the principal branch of the angle\footnote{
We use the following cartesian/polar notation for any \(w \in \mathbb C\): \(w = \Re w + i \Im w = \abs{w} \exp(i \Arg w)\), where the principal branch of argument is \(\Arg : \mathbb C \to [-\pi, \pi)\).  Additionally, the complex conjugate is denoted by \(\conj w\).
} \(\Arg \mathring{\mu}\) is a real-valued \emph{bounded and positive} function \(\Arg \mathring{\mu}(z) : \mathbb{D} \to [0,\pi)\).

The class of complex-valued analytic functions of \(\mathbb D\) with a positive real part is known as the \emph{Carath\'eodory class}. When \(F(z)\) is Carath\'eodory and its real part is furthermore bounded by a positive constant \(F^{+}\), the Riesz--Herglotz theorem (\cite{budisic_conditioning_2012},~\cite[\S 12.10]{henrici_applied_1977}) asserts that \(F(z)\) has an integral representation through measure \(\phi(\theta)\) on the circle \(\partial \mathbb D\):
\begin{equation}
  \label{eq:riesz-herglotz}
  F(z) = i \Im F(0) + \frac{1}{2\pi} \int_{-\pi}^{\pi} \frac{e^{i\theta} + z}{e^{i\theta}-z}d\phi(\theta).
\end{equation}
The integral in~\eqref{eq:riesz-herglotz} is very closely related to the analytic representation~\eqref{eq:analytic-representation} of the measure \(\phi\) as it can be rewritten to give
\begin{equation}
  \frac{1}{2\pi}\int_{-\pi}^{\pi} \frac{e^{i\theta} + z}{e^{i\theta}-z}d\phi(\theta) = -\tau_{\phi}(0) - 2i \mathring{\phi}(z).
\end{equation}
We can, therefore, re-state the representation formula for \(F(z)\) as
\begin{equation}
  \label{eq:riesz-herglotz-analytic}
  F(z) = i \Im F(0) - \tau_{\phi}(0) - 2i \mathring{\phi}(z).
\end{equation}
Riesz--Herglotz theorem additionally establishes that the measure \(\phi\) is continuous, i.e., there exists an integrable density function \(\phi'(\theta)\). Furthermore, the bounds on the real part of \(F(z)\),  \(\Re F(z) \in [0, F^{+})\), also pointwise bound the density, implying \(\phi'(\theta) \in [0, F^{+})\).

As shown earlier, \(\Arg \mathring{\mu}(z) = -i \Log \mathring{\mu}(z)\) is a positive real-valued function bounded by \([0,\pi)\),  Riesz--Herglotz machinery can be applied to the principal branch of logarithm of \(\mu(z)\). For any \(M \geq 0\) (to be chosen later) set
\begin{equation}\label{eq:choice-of-F}
F(z) := -i\Log\left[ \mathring{\mu}(z) + i M \right].
\end{equation}

At \(z = 0\) we obtain
\begin{equation}\label{eq:F-at-0}
  \begin{aligned}[t]
    F(0) &= -i \Log[\mathring{\mu}(0) + iM] = -i \Log[i(\tau_{\mu}(0) + M)] \\
    &= \pi/2 - i \log[\tau_{\mu}(0) + M],
  \end{aligned}
\end{equation}
which implies that the Riesz--Herglotz representation is
\begin{equation}\label{eq:RH-given-F}
  F(z) = -i\log[\tau_{\mu}(0) + M] - \tau_{\phi}(0) - 2i\mathring{\phi}(z).
\end{equation}
Equating~\eqref{eq:choice-of-F} and~\eqref{eq:RH-given-F} at \(z=0\), we find that the zeroth moment of \(\phi\) is always \(\tau_{\phi}(0) = \pi/2,\) for any measure \(\mu\).

From equivalence of~\eqref{eq:choice-of-F} and~\eqref{eq:RH-given-F} we can solve for \(\mathring{\mu}(z)\) to obtain its \emph{analytic phase representation} in terms of \(\mathring{\phi}(z)\):
\begin{equation}\label{eq:analytic-phase-representation}
\mathring{ \mu }(z) = -iM  -i [\tau_{\mu}(0) + M] \exp[2 \mathring{\phi} (z)],
\end{equation}
for any choice of \(M \geq 0\).

The representation formula~\eqref{eq:analytic-phase-representation}
establishes a connection between an arbitrary measure \(\mu\), which can be
singular, and a measure \(\phi\) with a bounded density function \(\phi'\). As
\(\mathring{\phi}\) arises as an argument of a complex exponential
in~\eqref{eq:analytic-phase-representation}, we refer to the measure \(\phi\)
and density \(\phi'\) as the \emph{phase measure} and \emph{phase density},
respectively, of the measure \(\mu\).

Relations between three representations
(Figure~\ref{fig:three-representations}) apply to the phase measure \(\phi\)
as well, with the important distinction that the density \(\phi'\) exists and
is bounded, which brings us a step closer to applicability of
Plemelj--Sokhotski formulas for reconstruction. While continuity of \(\phi'\)
cannot be theoretically expected, its assured boundedness makes it possible to
approximate \(\phi'\) by a continuous function, to which Plemelj--Sokhotski
formulas then naturally apply.

\begin{figure}[htb]
  \centering
  \begin{tikzpicture}[
    description/.style={fill=white, inner sep=2pt},
    mathnode/.style={execute at begin node=\(\scriptstyle, execute at end node=\)},
    textnode/.style={execute at begin node=\footnotesize},
    ge/.style={fontscale=1.5}]
\matrix (m)  [matrix of math nodes, row sep=4em,
column sep=6em, minimum width=2em, cells=ge, anchor=center] 
{ \mu &  &  & \phi\\
      &  \mathring{\mu}(z) &  \mathring{\phi}(z) & \\
  \tau_{\mu}(k) &  & & \tau_{\phi}(k)\\};
\path[-]
(m-1-1) edge [dashed,thick,->] node [left,mathnode] {\int e^{-ik\theta}d\mu}
(m-3-1)
(m-2-2) edge [dashed,thick,->] node [sloped, anchor=center,below,mathnode] {\displaystyle \ilim_{z \to e^{i\theta}} - \elim_{z \to e^{i\theta}}} (m-1-1)
(m-2-2) edge [dashed,thick,] node [sloped, anchor=center,above,textnode] {Plemelj--Sokhotski} (m-1-1)
(m-2-3) edge [dashed,thick,->] node [right, sloped, anchor=center,above,textnode] {Power expansion}
(m-3-4)
(m-2-2) edge [dashed,thick,->,transform canvas={yshift=-.5ex}] node [below,mathnode] {\Log}
(m-2-3)
(m-2-3) edge [dashed,thick,->,transform canvas={yshift=.5ex}] node [above,mathnode] {\exp}
(m-2-2)
(m-3-4) edge [thick,->] node [right,textnode] {\parbox[c]{5em}{Closure: MAXENT}}
(m-1-4)
(m-3-1) edge [dashed,thick,->] node [sloped, anchor=center,above,textnode] {Formal series}
(m-2-2)
(m-1-4) edge [dashed,thick,->] node [sloped, anchor=center,above,textnode] {Cauchy Transform}
(m-2-3)
(m-1-4) edge [dashed,thick] node [dashed,sloped, anchor=center,below,mathnode] {\int (1-e^{-i\theta}z)^{-1}d\phi}
(m-2-3)
(m-3-1) edge [thick,->] node [anchor=center,below,mathnode] {}
(m-3-4)
(m-3-1) edge node [anchor=center,above, textnode] {Triangular moment transformation}
(m-3-4)
(m-1-4) edge [thick,->] node [anchor=center,above, textnode] {Hilbert transform of density \(\mathcal H \phi'\) }
(m-1-1)
(m-1-4) edge node [anchor=center,below, mathnode] {\exp(\mathcal H \phi') \sin(\phi')}
(m-1-1)
;
  \end{tikzpicture}
  \caption{Relations between measures \(\mu\) and \(\phi\), their moments, and their analytic representations. Full arrows represent steps actually performed in our implementation, while dashed arrows represent analytical justifications that are not numerically evaluated (cf.~\autoref{fig:three-representations}).}
  \label{fig:measure-phase-relations}
\end{figure}
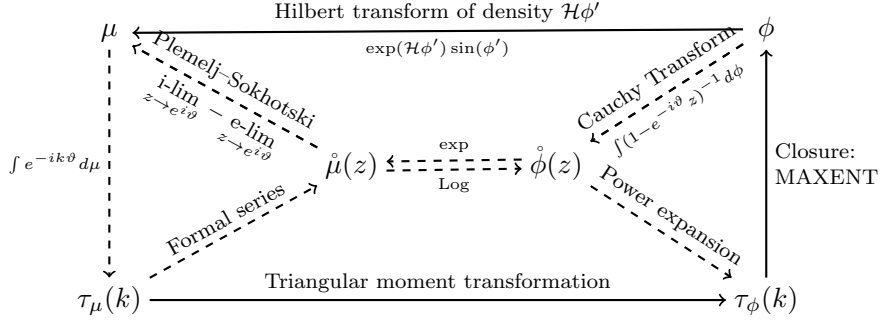

\autoref{fig:measure-phase-relations} illustrates the relations between discussed objects. In the next three sections, we operate in an applied context and describe how:
\begin{enumerate}[(a)]
\item phase moments \(\tau_{\phi}(k)\) are obtained from \(\tau_{\mu}(k)\) when \(k = 0,1,\dots,K \leq \infty\),
\item \(\phi'(\theta)\) is approximated by a continuous function from \(\tau_{\phi}(k)\) via entropy maximization, and
\item \(\mu\) is approximated using  \(\phi'(\mu)\) and its Hilbert transform \(\mathcal H \phi'(\mu)\).
\end{enumerate}

\subsection{A triangular transformation conditions the moment sequence}
\label{sec:conditioning}

As the expression~\eqref{eq:analytic-phase-representation} contains analytic
functions, it is obvious that the full sequence of phase moments
\(\tau_{\phi}(k), \forall k \in \mathbb Z\) can be computed from the sequence
\(\tau_{\mu}(k)\) \(\forall k \in \mathbb Z\). A subtler point is that to
compute the first \(K\) phase moments, \(\tau_{\phi}(k),\, 0 \leq k \leq
K-1\), we need just the first \(K\) moments \(\tau_{\mu}(k),\, 0 \leq k \leq
K-1\). The analytic phase
representation~\eqref{eq:analytic-phase-representation} connects power
expansions of \(\mathring{\mu}\) and \(\mathring{\phi}\) at \(z=0\)
\begin{align*}
  M + \sum_{n=0}^{\infty} \tau_{\mu}(n) z^{n} &= -[\tau_{\mu}(0) +M] \exp\left[ 2i \sum_{k=0}^{\infty} \tau_{\phi}(k) z^{k} \right]. \\
\Log\left[M + \sum_{n=0}^{\infty} \tau_{\mu}(n) z^{n}\right] &= \pm i\pi + \log[\tau_{\mu}(0) + M] + 2i \sum_{k=0}^{\infty} \tau_{\phi}(k) z^{k}.
\end{align*}
The choice of the sign is made by evaluation at \(z=0\) with \(\tau_{\phi}(0) = \pi/2\), resulting in
\begin{equation}\label{eq:moment-transform}
\sum_{k=0}^{\infty}\tau_{\phi}(k) z^{k} = \frac{\pi}{2} + \frac{i}{2} \log [\tau_{\mu}(0) + M] - \frac{i}{2} \Log\left[M + \sum_{n=0}^{\infty}\tau_{\mu}(n)z^{n}\right].
\end{equation}

For the choice of \(M = 1\), \(\Log\) can be directly expanded into formal Mercator series \(\Log (1 + z) = -\sum_{n=1}^{\infty} \frac{(-1)^{n}}{n}z^{n}\):
\begin{equation}\label{eq:moment-transform-m-1}
\sum_{k=0}^{\infty}\tau_{\phi}(k) z^{k} = \frac{\pi}{2} + \frac{i}{2} \log [\tau_{\mu}(0) + 1] + \frac{i}{2} \sum_{k=1}^{\infty} \frac{(-1)^{k}}{k}\left[ \sum_{n=0}^{\infty} \tau_{\mu}(n) z^{n}\right]^{k}.
\end{equation}
Alternatively, we can extract \(M + \tau_{\mu}(0)\) from the argument of \(\Log\) in~\eqref{eq:moment-transform} to obtain:
\begin{equation}\label{eq:moment-transform-m-0}
  \begin{aligned}
    \sum_{k=0}^{\infty}\tau_{\phi}(k) z^{k} &= \frac{\pi}{2} - \frac{i}{2} \Log\left[1 + \sum_{n=1}^{\infty} \frac{\tau_{\mu}(n)}{M + \tau_{\mu}(0)}z^{n}\right] \\
 &= \frac{\pi}{2} + \frac{i}{2} \sum_{k=1}^{\infty} \frac{(-1)^{k}}{k} \left[ \sum_{n=1}^{\infty}\frac{\tau_{\mu}(n)}{M + \tau_{\mu}(0)}z^{n}\right]^{k}
  \end{aligned}
\end{equation}
While expressions~\eqref{eq:moment-transform-m-1} and~\eqref{eq:moment-transform-m-0} are largely similar, we draw the attention to the lower bound of the inner sum \(\sum_{n}\), which runs from \(n=0\) in~\eqref{eq:moment-transform-m-1} and from \(n=1\) in~\eqref{eq:moment-transform-m-0}. In the remainder of the paper, we will discuss only \(M = 0\) case of formulas, corresponding to~\eqref{eq:moment-transform-m-0}.

In either~\eqref{eq:moment-transform-m-1} or~\eqref{eq:moment-transform-m-0}
the exponents of \(z\) are non-negative, so coefficient \(\tau_{\phi}(k)\)
cannot depend on any coefficients \(\tau_{\mu}(n)\) in which \(n>k\). It
follows that to compute a finite number of \(\tau_{\phi}(k)\), only a finite
number of \(\tau_{\mu}(n)\) is needed.  Transformations of sequences \(\mathbf
u \mapsto \mathbf v\) in which element \(v_{n}\) depends only on elements
\(u_{1},u_{2},\dots,u_{n}\) are referred to as ``triangular'' in analogy to
linear transformations of this kind, which are represented by lower-triangular
matrices
\begin{equation*}
  \left[\begin{smallmatrix}
    v_{1} \\ v_{2} \\ \vdots \\ v_{n}
  \end{smallmatrix}\right] =
  \left[\begin{smallmatrix}
    \ast &  & & \\
    \ast &  \ast &  &\\
    \vdots &   & \ddots & \\
    \ast & \ast & \dots & \ast \\
  \end{smallmatrix}\right]
  \left[\begin{smallmatrix}
    u_{1} \\ u_{2} \\ \vdots \\ u_{n}
  \end{smallmatrix}\right].
\end{equation*}
Triangular transformations of sequences are important in calculations where only a finite number of input elements are available, as it implies that at least the initial part of the output sequence can be obtained exactly.

\subsection{A continuous approximation of the phase density using entropy maximization}
\label{sec:maxent}

Assume that a finite set of complex numbers \(\tau(k),\, k = 0, \dots, K-1\) is chosen such that there exists at least one measure \(\rho'(\theta)d\theta\) whose truncated moment sequence matches numbers \(\tau(k)\):
\begin{equation}
\tau_{\rho}(k) =  \tau(k),\quad k=0,\dots,K-1.\label{eq:truncated-moment-condition}
\end{equation}
The set of all measures that satisfy the truncated moment condition~\eqref{eq:truncated-moment-condition} is called the \emph{feasible set} \(\mathcal{F}\). Whenever \(K < \infty\), the feasible set is uncountable and the inverse problem, i.e., selecting a single measure that corresponds to the truncated moment conditions, does not have a unique solution.

Entropy maximizations (MAXENT) are a broad class of closures which are well-understood for various domains of measures~\cites{junk_maximum_2000,borwein_maximum_2012}. We discuss the simplest case of Boltzmann-Shannon entropy maximization on the continuous domain: real-valued measures on a real, bounded interval. Since the measure is real-valued, truncated moment conditions~\eqref{eq:truncated-moment-condition} additionally define the negative-order moments of orders up to \(-(K-1)\) through relation \(\tau_{-k} = \conj \tau_{k}\). The entropy maximization solves the inverse problem by selecting the measure that lies in the feasible set and maximizes the Boltzmann-Shannon entropy functional:
\begin{equation}
  \label{eq:entropy-optimization}
  \argmax_{\rho \in \mathcal{F}} \int \log \rho'(\theta) d\rho(\theta).
\end{equation}
In the setup as simple as ours, the solution is easily found by the Lagrange multipliers technique, searching for the stationary point of the argument of the integral in the functional
\begin{equation}
  \label{eq:lagrange-functional}
  \int \left[\rho'(\theta) \log \rho'(\theta)  + \frac{1}{2\pi}\sum_{k = -K}^{K} \alpha_{k} \rho'(\theta) e^{ik\theta}\right] d\theta,
\end{equation}
where \(\alpha_{k}\) are the Lagrange multipliers. As the result, we obtain the density
\begin{equation}
\rho'(\theta) = \exp\left( \sum_{k=-K}^{K} \alpha_{k} e^{ik\theta}\right),\quad \alpha_{-k} = \conj{\alpha_{k}}\label{eq:entropy-ansatz}
\end{equation}
for the solution of the inverse problem (with some constants absorbed into \(\alpha_{k}\)). The procedure can be made quite general, as in~\cite{borwein_maximum_2012,junk_maximum_2000}; here, we work with~\eqref{eq:entropy-ansatz} that states that solutions with maximal entropy are always exponential trigonometric polynomials whose degree matches the number of truncated moment conditions. Practically, this implies that the entropy maximization can be stated as a dual nonlinear optimization problem where coefficients \(\alpha_{k}\) are varied until the moments of measure \(\rho\) match the moment conditions. Applying entropy maximization defines higher moments \(\tau_{\rho}(k), \abs{k} > K\) through a recursive equation that depends on the coefficients \(\alpha_{k}\) as derived in~\cite{budisic_conditioning_2012}.

The entropy maximization is attractive for several reasons. First, in the
limit when a full moment sequence is known, MAXENT recovers the unique measure
that generates the moment sequence~\cite{junk_maximum_2000}. From the
practical standpoint, resolution of a truncated moment problem is attractive
in physical sciences since \emph{it always produces a smooth and positive
  density} through little more than a convex optimization
process. Unfortunately, not all finite sets of complex numbers \(\tau_{k}\)
result in truncated moment conditions that can be matched by the MAXENT
ansatz, as is easy to show by studying an atomic probability measure
\(\delta(\theta)\)~\cite{budisic_conditioning_2012}.

\subsection{The density is computed from the phase using the Hilbert transform}
\label{sec:reconstruction}

Using the entropy maximization as the closure of the truncated moment problem
for \(\phi\) allows us to assume that we have a smooth density
\(\phi'(\theta)\) featured in the analytic phase representation
~\eqref{eq:analytic-phase-representation}
\[
\mathring{\mu}(z) = -iM - i[\tau_{\mu}(0) + M] \exp\left[\frac{i}{\pi}\int_{-\pi}^{\pi} \frac{\phi'(\theta)d\theta}{1 - e^{-i\theta}z}\right].
\]
Plemelj--Sokhotski formulas~\eqref{eq:plemelj-sokhotski-formulas} can be used
to reconstruct the continuous density of \(\mu\) as a difference of
non-tangential limits of \(\mathring \mu(z)\). While our original measure
\(\mu\) did not necessarily have a continuous density, the use of entropy
maximization to obtain a smooth phase density has the effect of regularizing
\(\mu\) as well.

Due to analyticity of the exponential function, non-tangential limits of
\(\mathring{\mu}(z)\) can be stated in terms of limits of
\(\mathring{\phi}(z)\)
\begin{equation}
  \label{eq:non-tangential-of-mu}
  \begin{aligned}
    \ilim_{z \to e^{i\theta}}\mathring{\mu}(z) &= -iM - i[\tau_{\mu}(0)
    + M] \exp\left[ 2 \ilim_{z \to
        e^{i\theta}}\mathring{\phi}(z)\right], \\
    \elim_{z \to e^{i\theta}}\mathring{\mu}(z) &=
    i\tau_{\mu}(0) + \conj{\ilim_{z \to e^{i\theta}}\mathring{\mu}(z)} \\
    &= i\tau_{\mu}(0) + iM + i[\tau_{\mu}(0)
    + M] \exp\left[ 2 \conj{\ilim_{z \to
        e^{i\theta}}\mathring{\phi}(z)}\right] \\
    &= i\tau_{\mu}(0) + iM - i[\tau_{\mu}(0)
    + M] \exp\left[ 2\elim_{z \to
        e^{i\theta}}\mathring{\phi(z)}\right],
  \end{aligned}
\end{equation}
where we used~\eqref{eq:elim-as-ilim} and \(\tau_{\phi}(0) = \pi/2\) to relate the interior and exterior non-tangential limits.

Taking the difference of non-tangential limits for \(\mathring{\mu}(z)\) and
factoring the exponential yields \emph{the inversion formula}
\begin{align}
  \mu'(\theta) &= -2M- \tau_{\mu}(0) + 2[\tau_{\mu}(0) + M] e^{ -\mathcal{H}\phi'(\theta)} \, \sin \phi'(\theta).\label{eq:inversion-formula}
\end{align}
Note that the choice of \(M \geq 0\) implicitly affects moments
\(\tau_{\phi}(k)\), through power expansion
relation~\eqref{eq:moment-transform}, and, consequently, the recovered
\(\phi'(\theta)\).

While the inversion formula ~\eqref{eq:inversion-formula} seemingly involves
another singular integral through the Hilbert transform, this transformation
has a simple spectral representation. If a general function \(f\) is expanded
into Fourier series as \(f(\theta) = \sum_{n = -\infty}^{\infty} \tau_{n}
e^{in\theta}\), then the expansion of its Hilbert transform \(\mathcal Hf\) is
\begin{equation}
\label{eq:hilbert-fourier}
\mathcal{H}f(\theta) = -i\sum_{n = -\infty}^{\infty} \sgn(n)\tau_{n} e^{in\theta}.
\end{equation}
In other words, the Hilbert transform can be efficiently computed using
commonly available Fast Fourier Transform (FFT) and inverse FFT
algorithms.\footnote{Alternatively, if Chebyshev expansion of \(f(\theta)\) is
  given \(f(\theta) = \sum_{n} c_{n}T_{n}(\theta)\), then the Hilbert
  transform of a rescaled function \(f(\theta)(1-\theta^{2})^{-1/2}\) has an
  expansion \( \mathcal{H}\left[ f(\theta)(1-\theta^{2})^{-1/2} \right] =
  -\sum_{n} c_{n}U_{n-1}(\theta), \) where \(T_{n}\) and \(U_{n}\) are
  Chebyshev polynomials of the first and second kind, respectively. Additional
  formulas for an explicit expansion of \(\mathcal H f\) into Chebyshev
  polynomials of the first kind are derived over several different domains
  in~\cite{olver_computing_2011}.}

\subsection{An essential example: a point mass measure}
\label{sec:delta}

Although the high-level concepts of the moment conditioning are easily understood from \autoref{fig:measure-phase-relations}, the detailed calculations require delicate choices of branches of inverse functions, which can be performed explicitly on an atomic measure \(\mu\). When \(\mu(\theta)\) is a point-mass at \(\theta = a\), i.e., \(d\mu(\theta) = \delta_{a}(\theta)d\theta\), most of the expressions can be determined analytically as integrals against the Dirac-\(\delta_{a}\), which simplifies to evaluations of the integrand at \(\theta = a.\)

Trigonometric moments are given by
\(  \tau_{\mu}(k) = (2\pi)^{-1}\int_{-\pi}^{\pi} e^{-ik\theta}\delta_{a}(\theta)d\theta = e^{-ika}/(2\pi)\) and the
analytic representation \(\mathring{\mu}(z)\) on \(\mathbb D\) is
\begin{equation}
\mathring{\mu}(z) = \frac{i}{2\pi}\int_{-\pi}^{\pi} \frac{\delta_{a}(\theta)d\theta}{1-e^{-i\theta}z} = \frac{i}{2\pi} \frac{1}{1 - e^{-ia}z}
= i \sum_{n=0}^{\infty} \frac{e^{-ian}}{2\pi} z^{n}.
\end{equation}

To calculate \(\mathring{\phi}(z)\), use the analytic phase representation ~\eqref{eq:analytic-phase-representation} and write
\begin{align*}
  \exp[2 \mathring{\phi} (z)] &= \frac{i \mathring{ \mu }(z) - M}{\tau_{\mu}(0) + M}
  = \frac{2\pi i}{1 + 2\pi M} \mathring{\mu}(z) - \frac{2\pi M}{1 + 2\pi M} \\
  &= \frac{1}{1 + 2\pi M} \cdot \frac{1 +2\pi M - 2\pi M e^{-ia}z}{e^{-ia}z - 1}
\end{align*}

At this point, choose \(M = 0\) which reduces the representation to
\begin{align*}
  \exp[2 \mathring{\phi} (z)] &= \frac{1}{e^{-ia}z - 1}.\\
  \mathring{\phi} (z) &= -\frac{1}{2}\Log[e^{-ia}z - 1].
\end{align*}
where the branch \(\Log[-1] = -i\pi\) is preferred so that \(\mathring{\phi}(z)\) would match the value of the zeroth moment \(\tau_{\phi}(0) = \pi/2\) as calculated in~\eqref{eq:RH-given-F}.

As \(\mathring{\phi}(z)\) is analytic inside \(\mathbb D\), calculating \(\ilim_{z \to e^{i\theta}}\) is done trivially by substituting the limit point. Taking \(\elim\) requires using the associate function, as in~\eqref{eq:elim-as-ilim}, to obtain
\begin{equation}
  \label{eq:nontangential-delta}
  \begin{aligned}
    \ilim_{z \to e^{i\theta}} \mathring{\phi}(z) &=
-\frac{1}{2}\Log[e^{i(\theta-a)} - 1] \\
    \elim_{z \to e^{i\theta}} \mathring{\phi}(z) &=
-\frac{1}{2}\Log[e^{-i(\theta-a)} - 1] + i\frac{\pi}{2}.
  \end{aligned}
\end{equation}

By Plemelj--Sokhotski~\eqref{eq:plemelj-sokhotski-density}, we obtain
\begin{equation}
  \label{eq:delta-density}
  \begin{aligned}
    \phi'(\theta) &= -i \left[\ilim_{z \to e^{i\theta}}
      \mathring{\phi}(z) - \elim_{z \to e^{i\theta}}
      \mathring{\phi}(z)\right] = \frac{i}{2} \Log \frac{e^{i(\theta-a)} - 1}{e^{-i(\theta-a)} - 1} - \frac{\pi}{2}\\
    &= -\Arccot \cot \frac{\theta-a}{2} + \pi
    = -\bmod\left( \frac{\theta-a}{2}, \pi \right) + \pi \\
    \mathcal{H} \phi'(\theta) &= -\left[\ilim_{z \to e^{i\theta}}
      \mathring{\phi}(z) + \elim_{z \to e^{i\theta}}
      \mathring{\phi}(z)\right] +i\tau_{\phi}(0) \\&= \frac{1}{2}\Log[2 - e^{i(\theta-a)} - e^{-i(\theta-a)}] = \Log\left[ 2 \abs{\sin\frac{\theta-a}{2}} \right]
  \end{aligned}.
\end{equation}
In all calculations above, care must be taken to select appropriate branches that ensure that \(\phi'(\theta) \in [0, \pi],\) as stipulated by the Riesz--Herglotz formula.

To evaluate \(\mu'(\theta)\), we use the inversion formula~\eqref{eq:inversion-formula}
\begin{align*}
    \mu'(\theta) &= - \tau_{\mu}(0) + 2\tau_{\mu}(0) e^{ -\mathcal{H}\phi'(\theta)} \, \sin \phi'(\theta) \\
    &= - \tau_{\mu}(0) + \tau_{\mu}(0) \frac{1}{\abs{\sin\frac{\theta-a}{2}}} \sin\left[ -\bmod\left( \frac{\theta-a}{2},\pi\right) + \pi\right] \\
    &= - \tau_{\mu}(0) + \tau_{\mu}(0)
\frac
{\sin \bmod\left( \frac{\theta-a}{2},\pi\right) } {\abs{\sin\frac{\theta-a}{2}}}.
\end{align*}
Note that for all \(\theta \not = 0\) it holds that \(\sin \bmod\left( \frac{\theta-a}{2},\pi\right) = \abs{\sin\frac{\theta-a}{2}},\) so \(\mu'(\theta \not = 0) \equiv 0.\) However, since \(\int d\mu(\theta) = 1\), it has to hold that \(\mu'(\theta) = \delta_{a}(\theta).\)

\section{Numerical examples}
\label{sec:numerics}

To demonstrate the effect that the conditioning has on solution of the inverse moment problem, we implemented the entire procedure as a MATLAB code, relying heavily on \texttt{Chebfun}~\cites{driscoll_Chebfun_2014,trefethen_approximation_2013} package.  \texttt{Chebfun} represents functions constructively using their Chebyshev or Fourier coefficients, instead of storing their pointwise values. Consequently, operations involving integrals and integral transforms are quicker and more precise than if it was implemented naively.

To condition the moments, we computed powers of the truncated analytic representation~\eqref{eq:moment-transform} directly, through discrete convolution of polynomial coefficients. Alternatively, it is possible to specify coefficients recursively using the Miller--Nakos algorithm ~\cite{nakos_expansions_1993}, which would become preferable for implementing a multivariate generalization of this method~\cite{budisic_conditioning_2012}.

The solution to the truncated moment problem using entropy maximization is cast as an unconstrained convex optimization problem and solved using MATLAB's internal derivative-free quasi-Newton optimization method (\texttt{fminunc})~\cite{the_mathworks_inc_matlab_2014}.  The optimization function is the square-error in first \(K\) moments  \(\mathcal{E} = \sum_{k=0}^{K-1} \abs{\tau(k) - \tau_{\rho}(k)}^{2}\), where \(\tau(k)\) represent input data, and \(\tau_{\rho}\) moments of the density in the maximum entropy ansatz~\eqref{eq:entropy-ansatz}. The optimization parameters are real and imaginary parts of complex coefficients \(\alpha_{k}\) for \(0 \leq k < K-1,\) where coefficients \(k < 0\) are obtained by imposing \(\conj{\alpha_{k}} = \alpha_{-k}\), as in~\eqref{eq:trig-moments-fourier-coeffs}.

Since \texttt{Chebfun} represents functions internally using their Fourier coefficients, implementing the Hilbert transform reduces to a multiplication by a vector of constants. To transform the coefficients, the relation~\eqref{eq:hilbert-fourier} is applied directly while \texttt{Chebfun} takes care of pointwise evaluation.

Numerical algorithm is applied to three different measures \(\mu\):
\begin{enumerate}[(a)]
\item a point mass measure,
  \begin{equation}
\mu'(\theta) = \delta_{a}(\theta),\label{eq:dirac-density}
\end{equation}
\item a continuous measure with smooth (superposition of Gaussians) density
  \begin{equation}
    \mu'(\theta) = 5 e^{ -(5\theta-10)^2 } + e^{ -(5\theta+7.5)^2 },
    \label{eq:smooth-density}
  \end{equation}
\item a continuous measure with discontinuous (rectangular) density
  \begin{equation}
    \mu'(\theta) =
    \begin{cases}
      1 & x \in \left(\frac{-\pi-1}{2}, \frac{\pi+1}{2}\right) \\
      0 & \text{otherwise}
    \end{cases}.
    \label{eq:disc-density}
  \end{equation}
\end{enumerate}
All the measures were further normalized to the unit mass.

Assume that only first \(K = 20\) initial moments \(\tau_{\mu}(k)\) of the ``true'' density \(\mu\) are available. Two approximations, \(\mu_{U}'\) and \(\mu_{C}'\), to the initial density are obtained using, respectively,  unconditioned MAXENT procedure and the conditioned MAXENT procedure (full arrows in Figure~\ref{fig:measure-phase-relations}). In the text, we will refer to them as U\-/MAXENT and C\-/MAXENT. Settings of the optimization algorithm that constructs maximum entropy density are kept the same in both procedure.

While theoretically U\-/MAXENT may not converge for singular measures, numerical implementations have a number of stopping criteria which prevent the code from running forever. In our case, MATLAB's stock optimization code terminated by reporting that the local minimum was numerically found in all but one example: when applying U\-/MAXENT to the singular measure, it stopped due to inability to reduce the moment error, without confirming numerically that the local minimum was found.

\begin{figure}[htb]
  \centering
  \begin{subfigure}[t]{0.3\linewidth}
    \centering
    \includegraphics[width=1.5in]{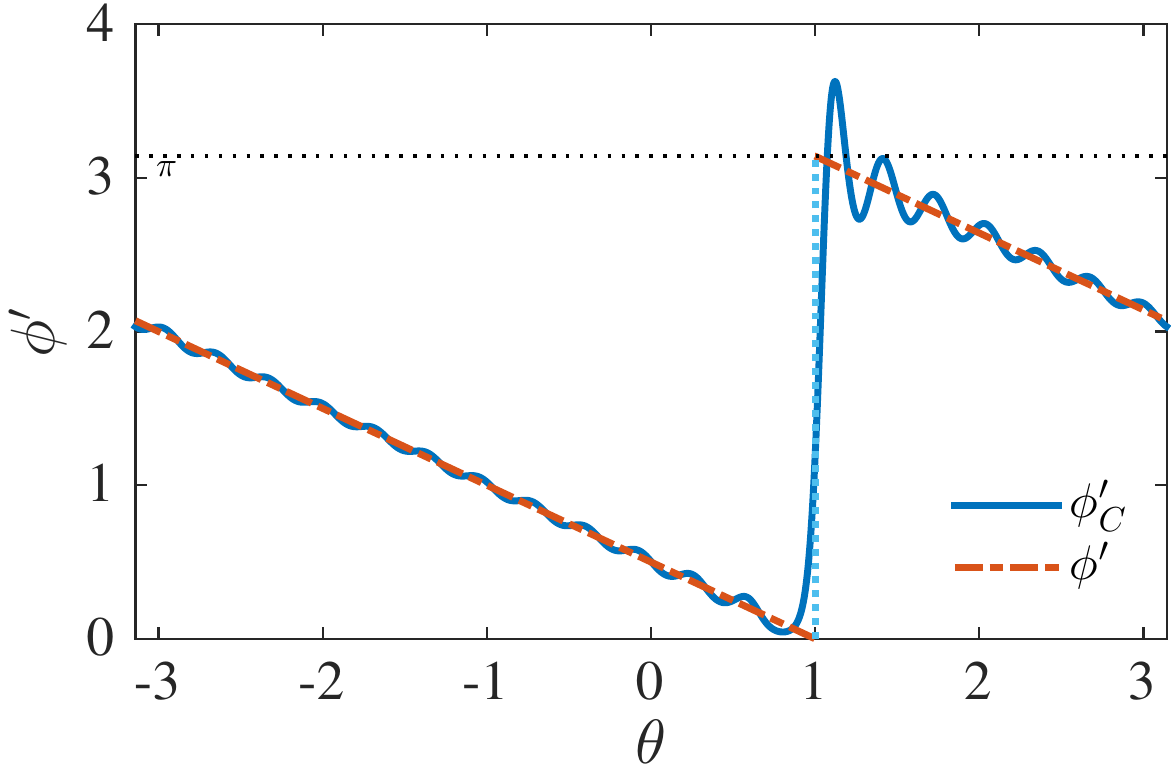}
    \caption{\(\phi'(\theta)\) and \(\phi_{C}'\) when \(\mu(\theta)\) is a point mass measure at \(1\).}
  \end{subfigure}\hfill
  \begin{subfigure}[t]{0.3\linewidth}
    \centering
    \includegraphics[width=1.5in]{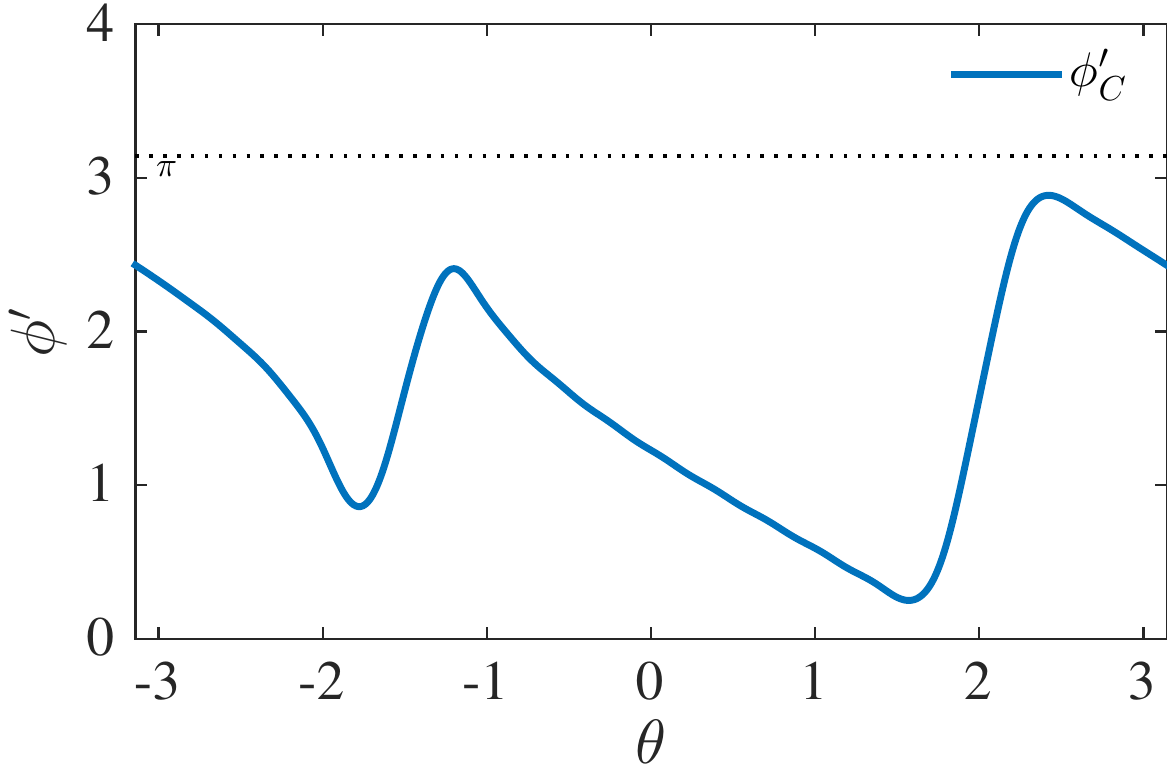}
    \caption{\(\phi'(\theta)\) when \(\mu'(\theta)\)  is a sum of Gaussians.}
  \end{subfigure}\hfill
  \begin{subfigure}[t]{0.3\linewidth}
    \centering
    \includegraphics[width=1.5in]{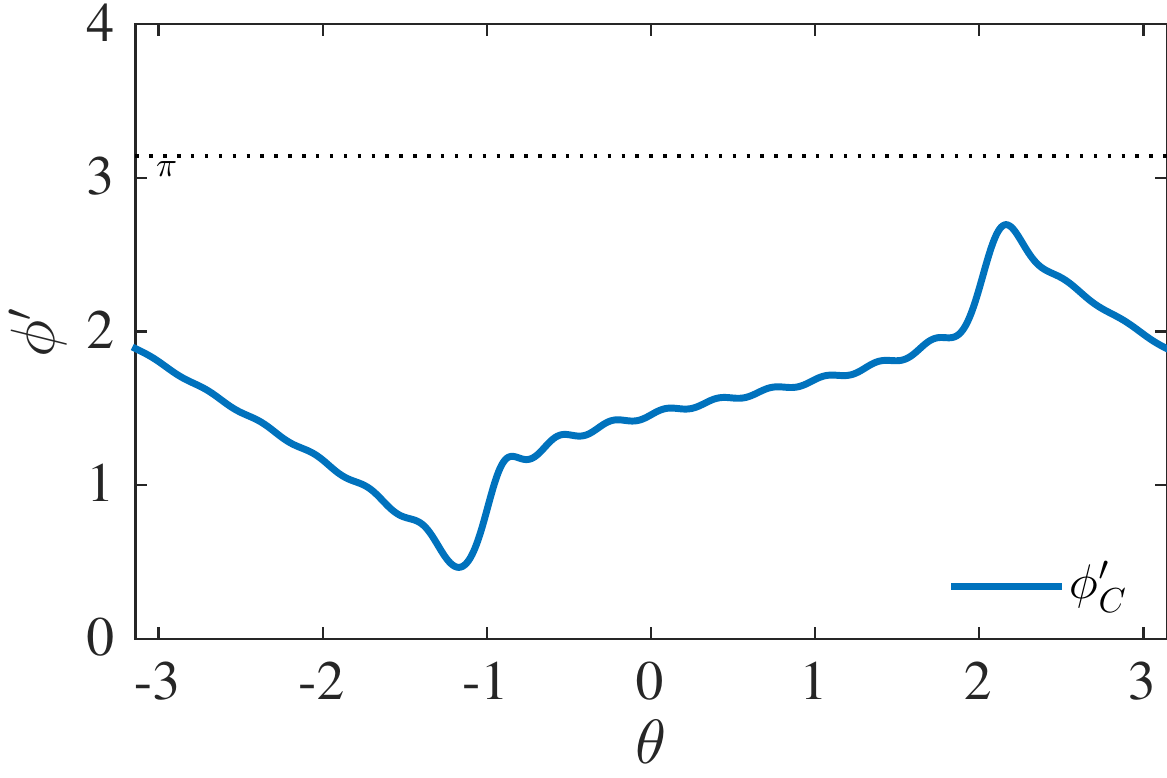}
    \caption{\(\phi'(\theta)\) when \(\mu'(\theta)\)  is rectangular.}
  \end{subfigure}
\caption{Reconstructed phase densities \(\phi_{C}'\) for three studied ``ground'' measures \(\mu\).}\label{fig:phase-densities}
\end{figure}

Before we discuss differences between U- and C\-/MAXENT, Figure~\ref{fig:phase-densities} shows what reconstructed phase densities \(\phi'\) look like for all three studied measures. While the phase density of \(\delta_{a}\) can be easily computed analytically (Section~\ref{sec:delta}), we have not attempted to do so for the remaining two densities. Loosely speaking, \(\phi'(\theta)\) appear to resemble ``integrals'' of \(\mu'(\theta)\) as singularities correspond to finite discontinuities, discontinuous changes are smoothed out, and constant regions are transformed into slopes. Intuitively, this should make \(\phi'\) smoother than \(\mu'\) and therefore more amenable to reconstruction by entropy maximization.

Since the reconstruction algorithm was performed in the Fourier space, it is not surprising that the Gibbs phenomenon~\cite{hewitt_gibbs-wilbraham_1979} at discontinuities results in overshooting of upper bound at \(\pi\), with effects on \(\mu_{C}\) discussed later in this section. Note that the maximum entropy ansatz \(e^{P(\theta)}\) produces overshooting only in the positive direction, because Fourier coefficients of exponent \(P(x)\) are the optimization parameters, and not the full sequence of Fourier coefficients of the density \(e^{P(\theta)}\).

\begin{figure}[htb]
  \centering
  \begin{subfigure}[t]{0.32\linewidth}
    \centering
    \includegraphics[width=1.5in]{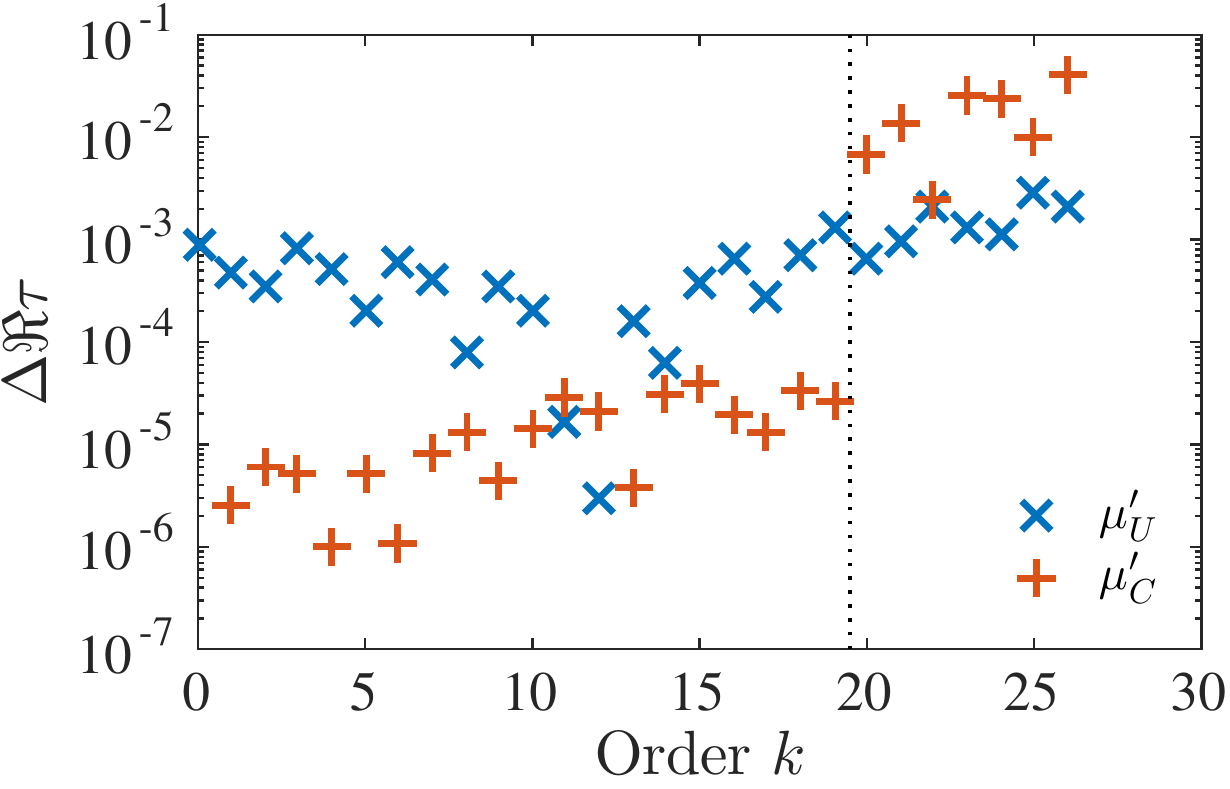}
    \caption{\(\Delta \Re \tau_{\mu}\) when \(\mu(\theta)\) is a point mass at \(1\).}
  \end{subfigure}\hfill
  \begin{subfigure}[t]{0.32\linewidth}
    \centering
    \includegraphics[width=1.5in]{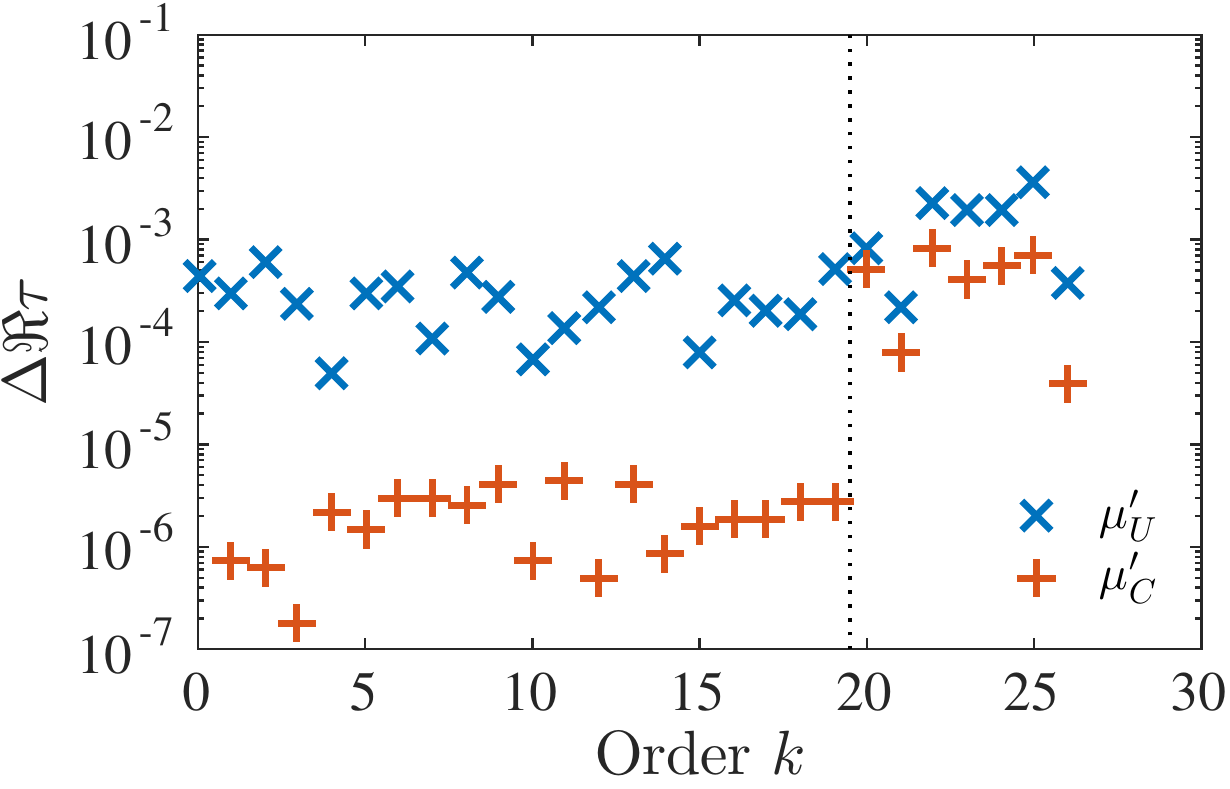}
    \caption{\(\Delta \Re \tau_{\mu}\) when \(\mu'(\theta)\)  is a sum of Gaussians.}
  \end{subfigure}\hfill
  \begin{subfigure}[t]{0.32\linewidth}
    \centering
    \includegraphics[width=1.5in]{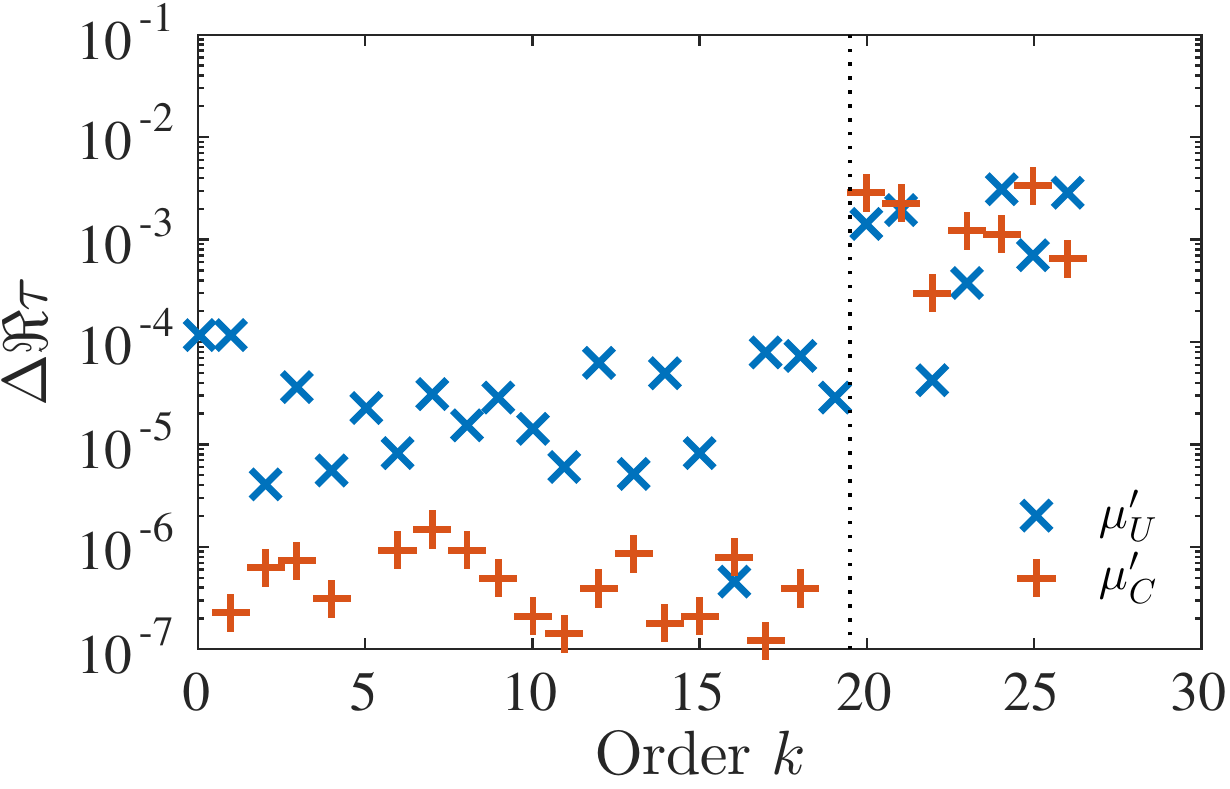}
    \caption{\(\Delta \Re \tau_{\mu}\) when \(\mu'(\theta)\)  is rectangular.}
  \end{subfigure}\\
  \begin{subfigure}[t]{0.32\linewidth}
    \centering
    \includegraphics[width=1.5in]{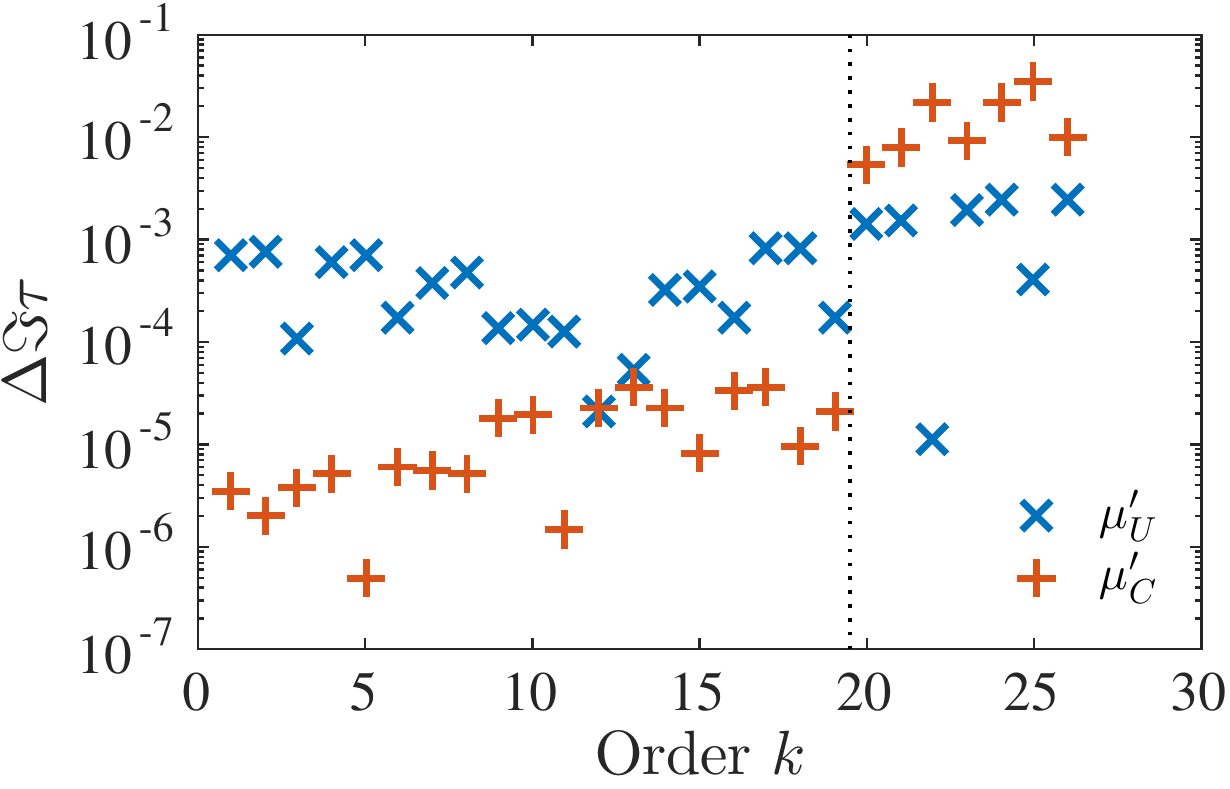}
    \caption{\(\Delta \Im \tau_{\mu}\) when \(\mu(\theta)\) is a point mass measure at \(1\).}
  \end{subfigure}\hfill
  \begin{subfigure}[t]{0.32\linewidth}
    \centering
    \includegraphics[width=1.5in]{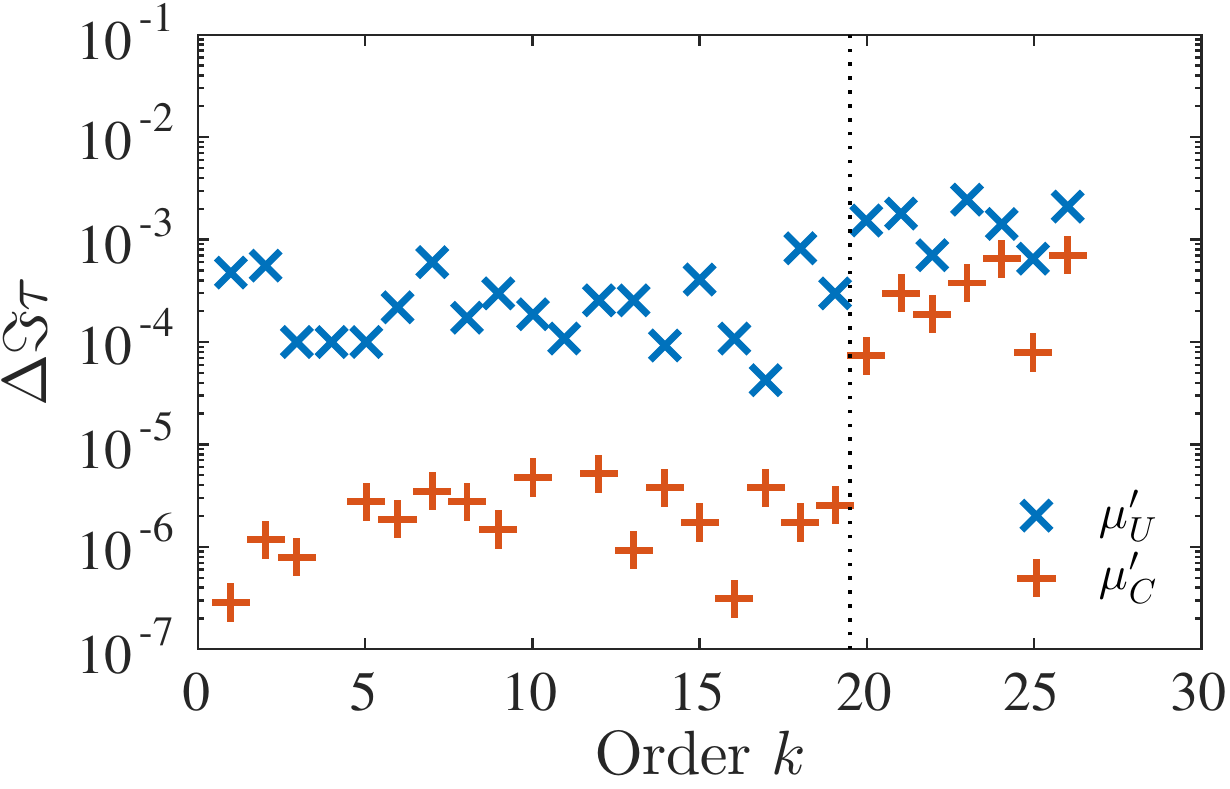}
    \caption{\(\Delta \Im \tau_{\mu}\) when \(\mu'(\theta)\)  is a sum of Gaussians.}
  \end{subfigure}\hfill
  \begin{subfigure}[t]{0.32\linewidth}
    \centering
    \includegraphics[width=1.5in]{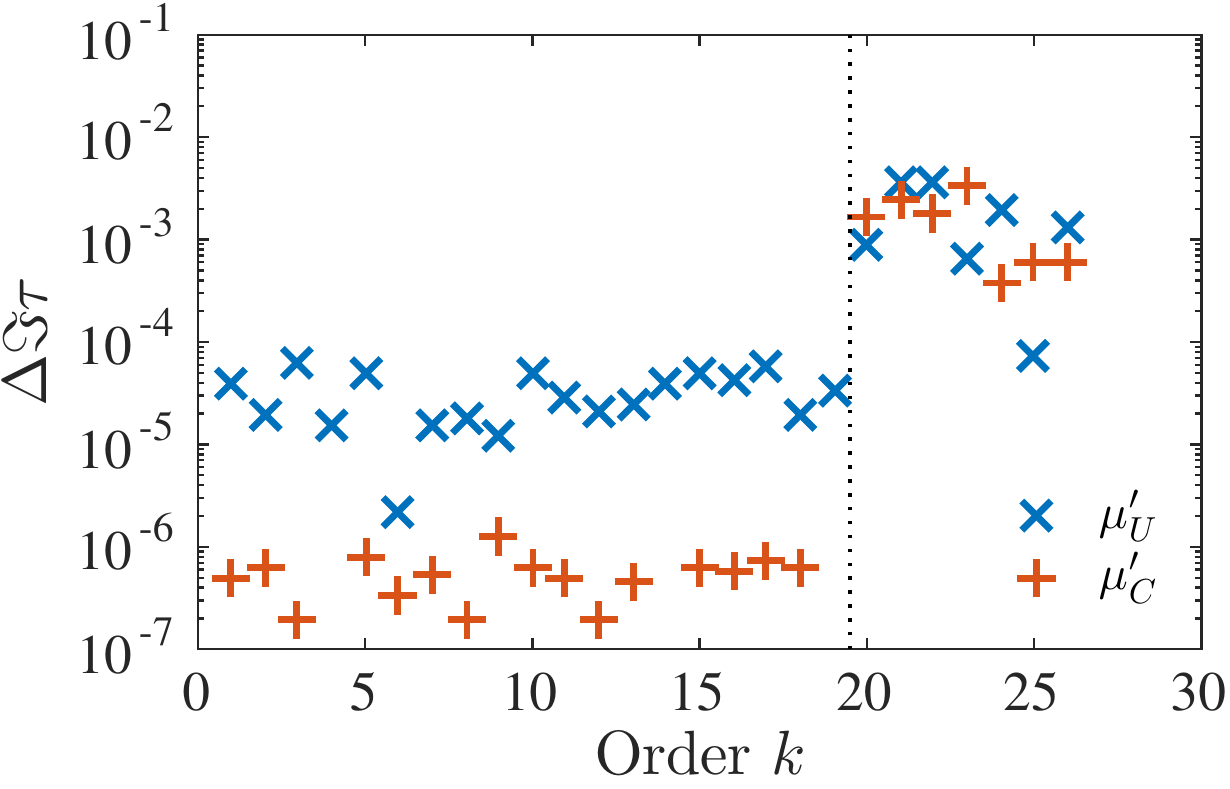}
    \caption{\(\Delta \Im \tau_{\mu}\) when \(\mu'(\theta)\)  is rectangular.}
  \end{subfigure}
  \caption{Errors of moments of measures reconstructed using MAXENT with and without conditioning, resp., \(\mu_{C}\) and \(\mu_{U}\), with respect to
    moments of the ``true'' measure \(\mu\) for the point mass measure~\eqref{eq:dirac-density}, a sum-of-Gaussians density~\eqref{eq:smooth-density}, and a discontinuous density~\eqref{eq:disc-density}.}~\label{fig:errors-in-moments}
\end{figure}

Figure~\ref{fig:errors-in-moments} compares moments of reconstructed densities, \(\mu_{C}'\) and \(\mu_{U}'\), with moments of \(\mu'\). Moment errors of \(\mu_{C}'\) at orders \(k < K = 20\) are lower by at least one order of magnitude than the same errors for \(\mu_{U}'\) in all studied examples. At orders \(k > K\), however, the advantage of \(\mu_{C}'\) is less clear.   Moment errors of \(\mu_{U}'\) remain roughly constant when \(k\) is increased beyond \(K\), except in the case of the discontinuous ``box'' density, where the error sharply increases. On the other hand, \(\mu_{C}'\) moment errors experience similar sharp increase in all three studied examples. Nevertheless, after the increase, the errors in moments of \(\mu_{C}'\) and \(\mu_{U}'\) are mostly comparable. The exception is the case of a point mass measure, where U\-/MAXENT outperforms C\-/MAXENT for \(k > K\).

\begin{figure}[htb]
  \centering
  \begin{subfigure}[t]{0.32\linewidth}
    \centering
    \includegraphics[width=1.5in]{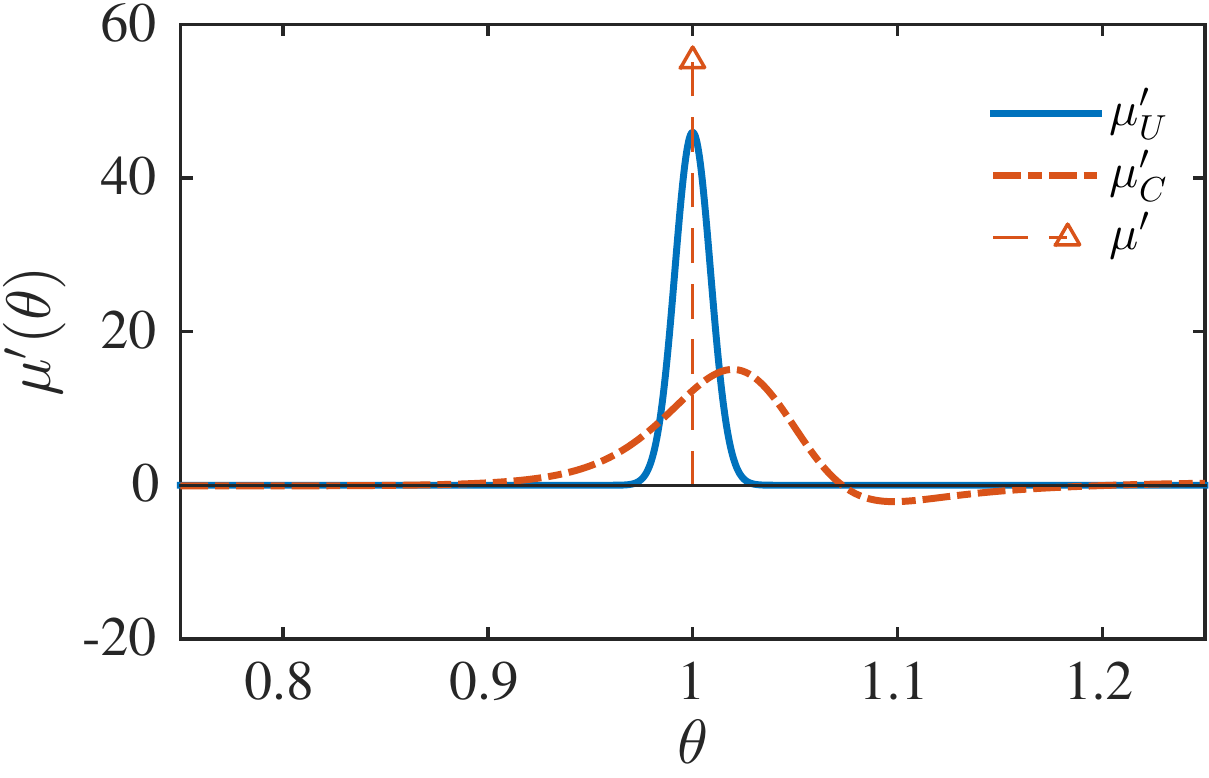}
    \caption{Pointwise density values when \(\mu(\theta)\) is a point mass measure at \(1\).}
  \end{subfigure} \hfill
  \begin{subfigure}[t]{0.32\linewidth}
    \centering
    \includegraphics[width=1.5in]{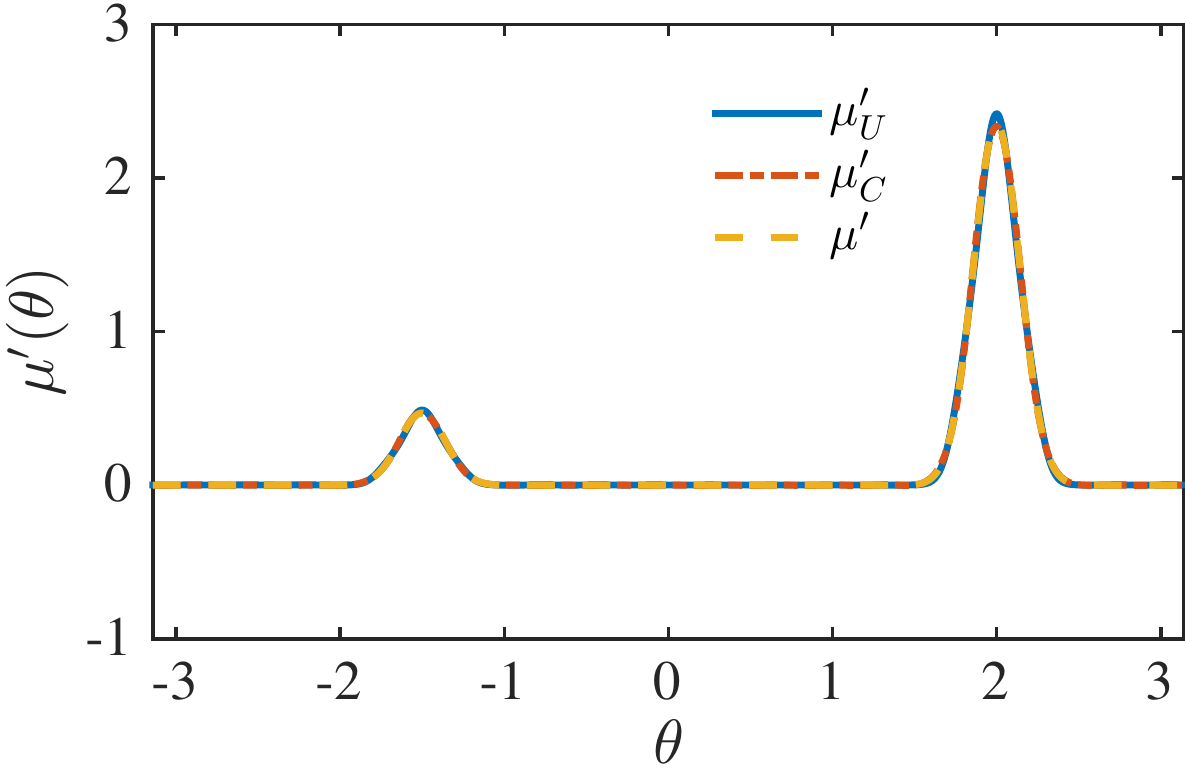}
    \caption{Pointwise density values when \(\mu'(\theta)\) is a sum of Gaussians.}
  \end{subfigure} \hfill
  \begin{subfigure}[t]{0.32\linewidth}
    \centering
    \includegraphics[width=1.5in]{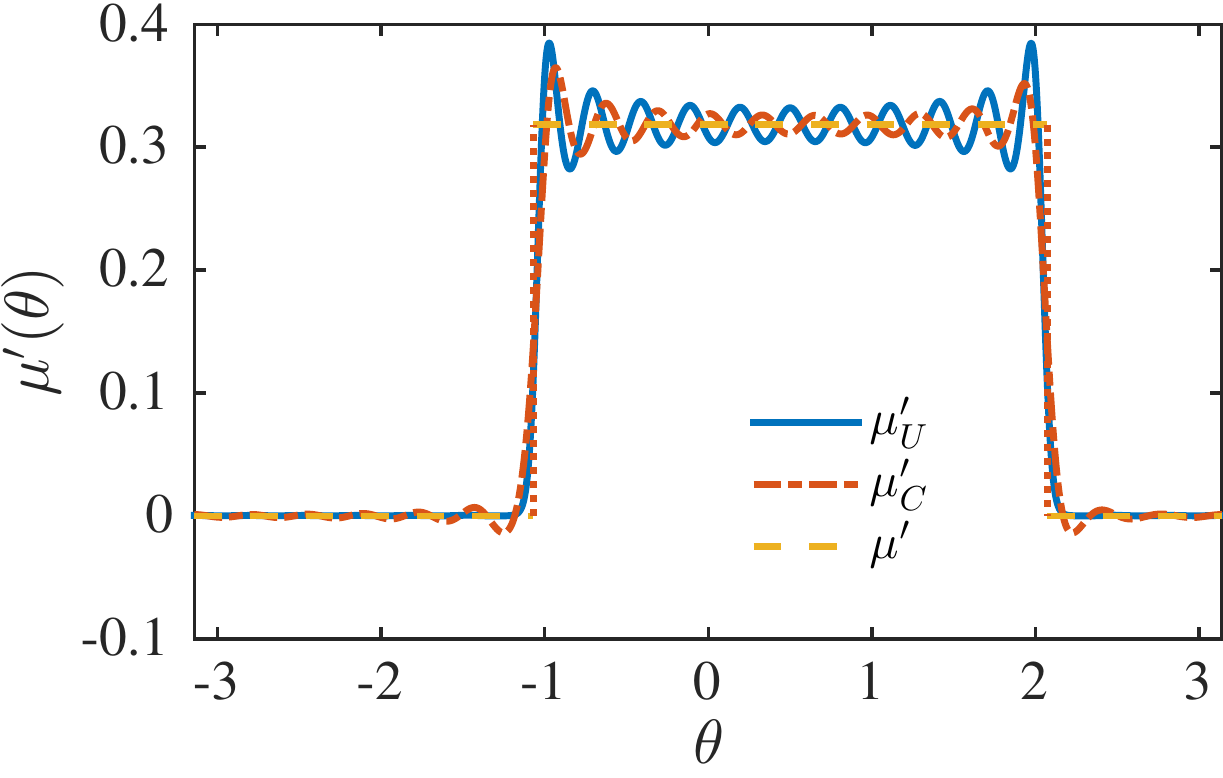}
    \caption{Pointwise density values when \(\mu'(\theta)\)  is rectangular.}
  \end{subfigure} \\
  \begin{subfigure}[t]{0.32\linewidth}
    \centering
    \includegraphics[width=1.5in]{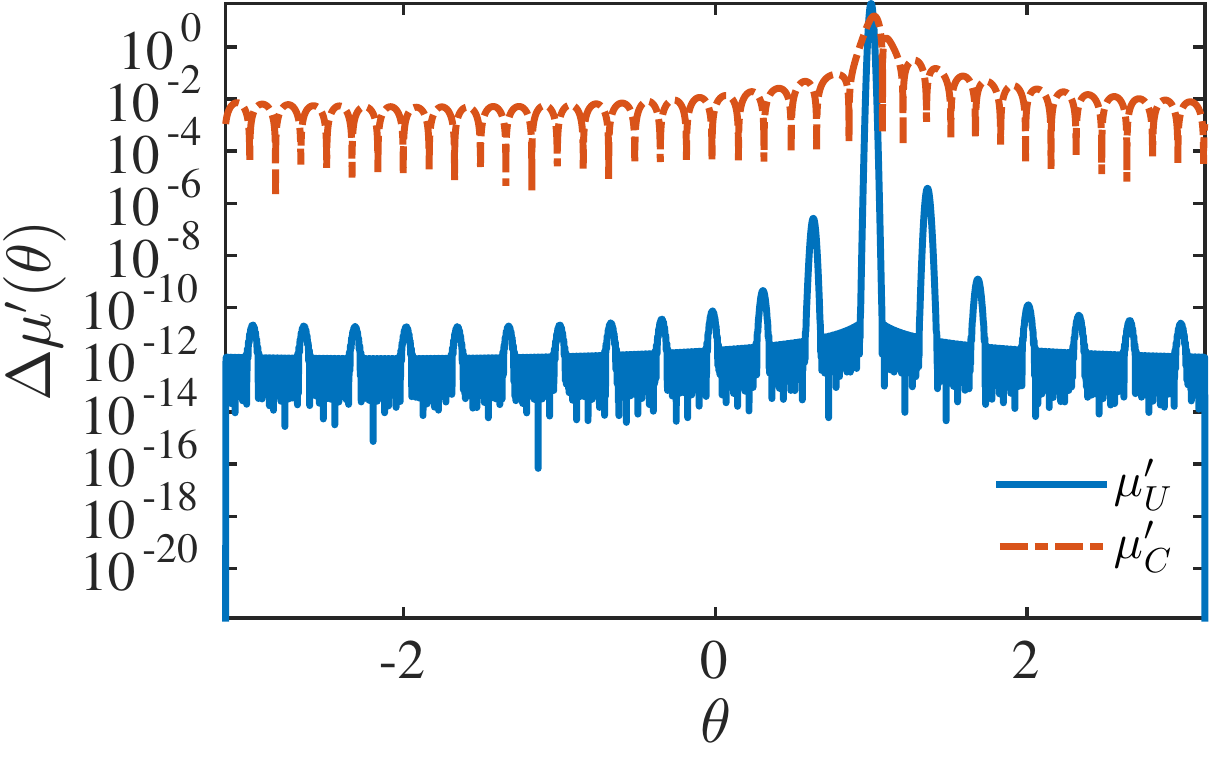}
    \caption{Pointwise density errors when \(\mu(\theta)\) is a point mass measure at \(1\).}
  \end{subfigure}\hfill
  \begin{subfigure}[t]{0.32\linewidth}
    \centering
    \includegraphics[width=1.5in]{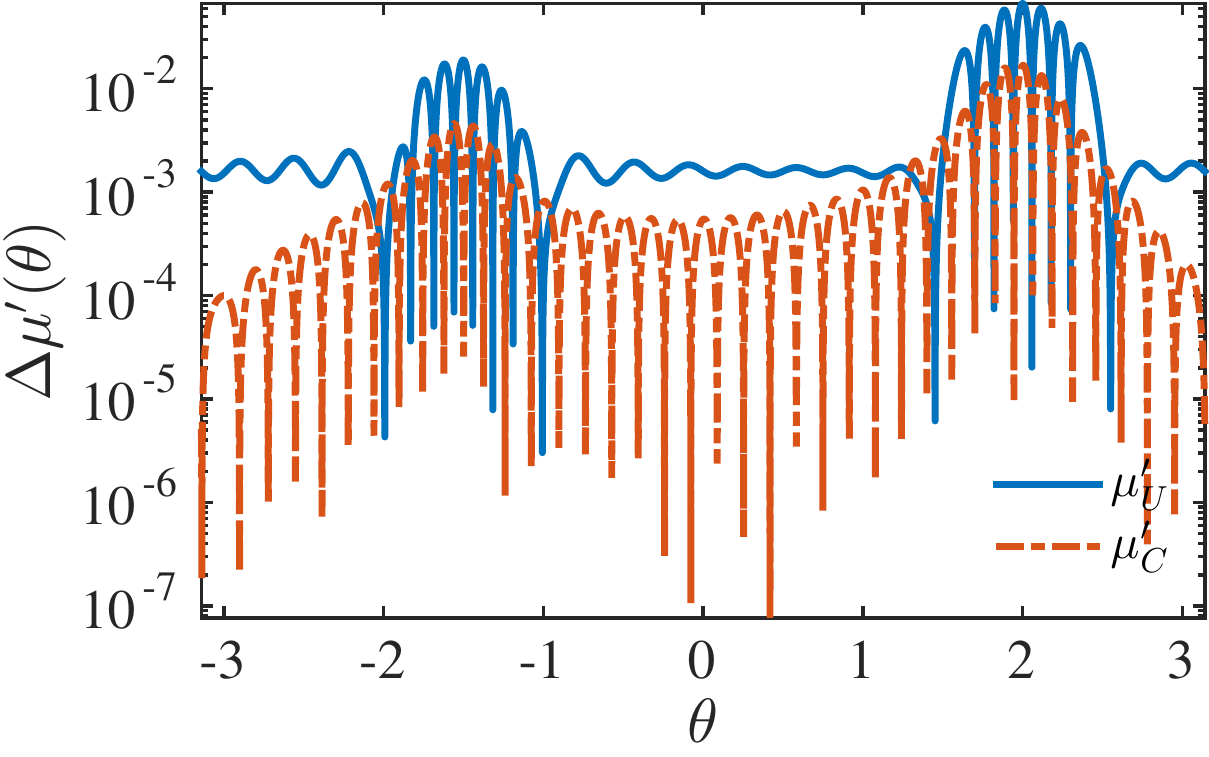}
    \caption{Pointwise density errors when \(\mu'(\theta)\) is a sum of Gaussians.}
  \end{subfigure}\hfill
  \begin{subfigure}[t]{0.32\linewidth}
    \centering
    \includegraphics[width=1.5in]{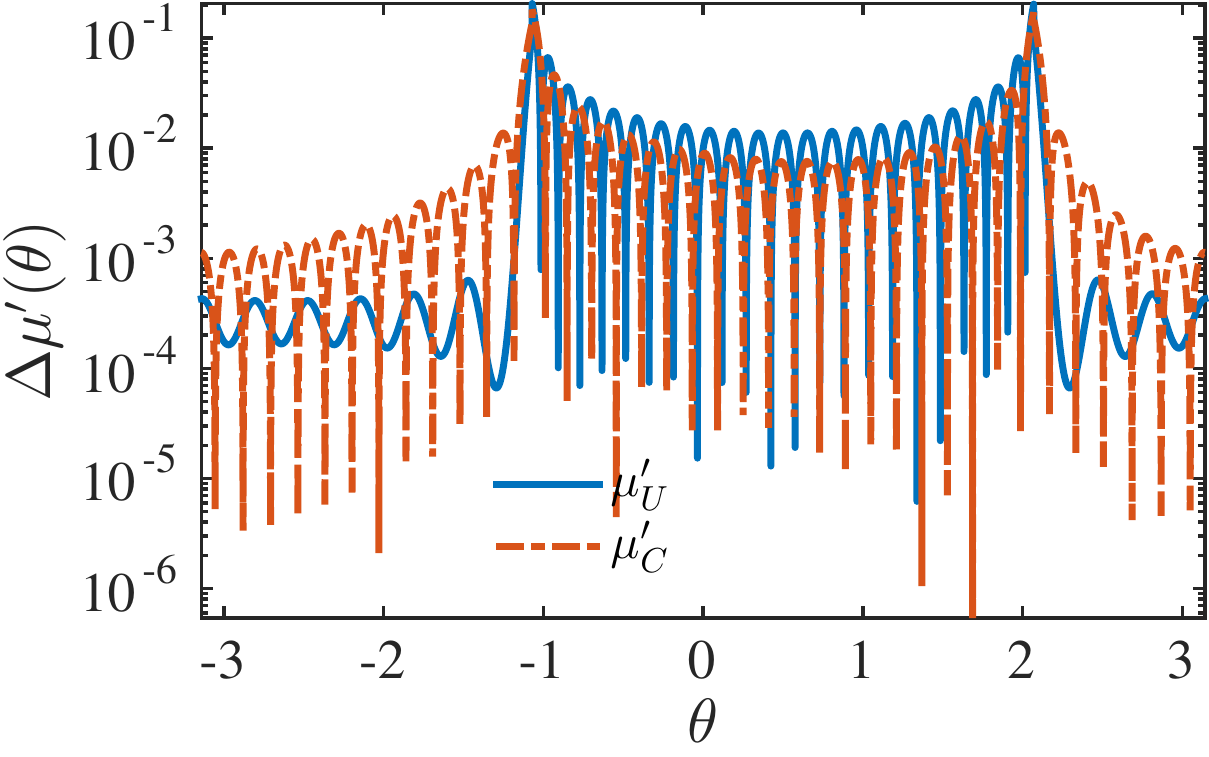}
    \caption{Pointwise density errors when \(\mu'(\theta)\)  is rectangular.}
  \end{subfigure}\hfill
\caption{Pointwise comparison of densities reconstructed using MAXENT with and without conditioning, resp., \(\mu'_{C}\) and \(\mu'_{U}\), with respect to
    the densities of ``true'' measures \(\mu\) for the point mass measure at \(1\)~\eqref{eq:dirac-density}, a sum-of-Gaussians density~\eqref{eq:smooth-density}, and a discontinuous density~\eqref{eq:disc-density}. \(\delta_{a}\) is indicated by an arrow, and for it, the error was not computed at \(\theta = a\) as it is always infinite.}\label{fig:pointwise-comparison}
\end{figure}

From pointwise comparison of densities in
Figure~\ref{fig:pointwise-comparison}, it can be seen that when \(\mu'\) is
smooth, like in the sum of Gaussians case, C\-/MAXENT provided a slightly
better approximation than U\-/MAXENT. It is not clear, though, whether that
difference is practically relevant.  On the other hand, the advantage to using
U\-/MAXENT, when the process converges, is that its density remains positive,
regardless of errors introduced. Reconstruction of densities in the case of
\(\delta_{a}\) and rectangular densities show that \(\mu'_{C}\) crosses \(0\)
into negative values. Both of these crossings are due to the Gibbs phenomenon,
which is manifested on two different levels. The undershoot of the point mass
reconstruction is a consequence of the Gibbs phenomenon in the reconstruction
of the phase \(\phi'\), as can be seen in Figure~\ref{fig:phase-densities},
which propagates through \(\sin \phi'(\theta)\) term in the inversion
formula~\eqref{eq:inversion-formula} and introduces negative values into
reconstruction \(\mu_{C}'\). On the other hand, the undershoot in the
reconstruction of the rectangular density does not have a counterpart in the
phase density, but is likely a direct Gibbs effect, as the inversion
formula~\eqref{eq:inversion-formula} was implicitly evaluated using finitely
many Fourier coefficients (see Table~\ref{tab:four-coeff}).

\begin{table}[htb]
  \centering
  \begin{tabular}[c]{l | c c c}
    Example & \(\mu_{U}'\) & \(\mu_{C}'\) & \(\phi_{C}'\) \\[1.5ex]
    Point mass & 1709 & 747 & 343 \\
    Sum of Gaussians & 291 & 253 & 193 \\
    Rectangular & 389 & 343 & 237 \\
  \end{tabular}
  \caption{Number of Fourier coefficients internally chosen by \texttt{Chebfun} toolbox to represent densities.}
  \label{tab:four-coeff}
\end{table}

\section{Conclusions}
\label{sec:conclusions}

Closures of inverse moment problems often require that the input moments
correspond to an absolutely continuous measure, as is the case of the maximum
entropy closure. If the input measure is singular, optimization algorithm
maximizing the entropy may not terminate. In this case, the finite moment sequence of the
measure can be conditioned so that the entropy maximization terminates and
resolves the inverse problem~\cite{budisic_conditioning_2012}. This paper
demonstrates that it is possible to numerically implement the conditioning
procedure for the case of measures on a compact interval and that, indeed,
conditioning of moments successfully allows the entropy maximization procedure
to terminate.

Numerical results show that conditioning moments before passing them to
entropy maximization consistently results in smaller errors between input
moments and the moments of the reconstruction, as compared to the
unconditioned entropy maximization. Comparing errors that an approximation
makes \emph{beyond} the first \(K\) moments shows less-definitive results:
conditioned and unconditioned procedures result in similar level of errors in
most cases, except when the measure is singular. In that case the
unconditioned algorithm surprisingly performed better than the conditioned
one. Pointwise comparison reveals that the final
step~\eqref{eq:inversion-formula} in the conditioned algorithm, inverting the
effect of moment conditioning, exposes the procedure to numerical errors. In
our implementation, this manifested as pointwise deterioration of the final
approximation, due to the Gibbs phenomenon. The measures chosen for the
examples were finitely determined by their moments, so there is no surprise
that even the unconditioned entropy maximization performed well.  In summary,
this paper successfully demonstrated that implementing the conditioning
procedure is practical, even though the used proof-of-concept algorithm may
not be the best possible implementation.

Future numerical analysis of our modification to the entropy maximization
would study the effects of error propagation through three major steps
(Figure~\ref{fig:measure-phase-relations}), as well as possibilities for
mitigating those errors. Additionally, extensions of the procedure into
higher-dimensional domains have been studied theoretically
in~\cite{budisic_conditioning_2012}, but the practical implementation would be
an important achievement. Extending conditioning of moments is almost
immediate, however the numerical implementation of the inversion step is less
straightforward; additional work is needed to see how this obstacle can be
overcome.

\bibliographystyle{spmpsci}      
\bibliography{moment-problem-ref}

\end{document}